\documentclass[12pt]{article}
\usepackage[letterpaper,margin=1in]{geometry}
\usepackage[utf8]{inputenc}
\usepackage[english]{babel}
\usepackage{enumitem}
\usepackage{ifthen}
\usepackage{subcaption}
\usepackage[]{lineno} 
\usepackage[table]{xcolor}
\usepackage[colorlinks=true,citecolor=blue,linkcolor=blue, urlcolor=gray]{hyperref}
\usepackage{fontawesome5}
\usepackage{xparse}       
\usepackage{etoolbox}     
\usepackage{authblk}      
\usepackage{xspace}
\usepackage{multirow}
\usepackage[font=small]{caption}

\makeatletter
\let\@authorline\@empty    
\let\@affilblock\@empty    
\newcommand{\emailicon}[1]{\href{mailto:#1}{\textsuperscript{\faEnvelope[regular]}}}
\newcommand{\orcidicon}[1]{\href{https://orcid.org/#1}{\textsuperscript{\faOrcid}}}
\newcommand{\homepage}[1]{\href{#1}{\textsuperscript{\faHome}}}
\newcommand{\authorentry}[5]{%
  \ifx\@authorline\@empty
  \else
    \g@addto@macro\@authorline{,\ }%
  \fi
  \g@addto@macro\@authorline{%
    #1%
    \ifx&#2&\else\textsuperscript{#2}\fi%
    \ifx&#3&\else\emailicon{#3}\fi%
    \ifx&#4&\else\orcidicon{#4}\fi%
    \ifx&#5&\else\homepage{#5}\fi%
  }%
}
\newcommand{\affilentry}[2]{%
  \g@addto@macro\@affilblock{%
    \textsuperscript{#1}#2\\%
  }%
}
\renewcommand{\maketitle}{%
  \begin{center}
    {\LARGE \bfseries \@title \par}
    \vskip 1em
    {\normalsize
      \@authorline
    }
    \vskip 1em
    {\small \@affilblock}
    \vskip 1em
    {\small \@date \par}
  \end{center}
  \vskip 2em
}
\makeatother

\usepackage{amsmath}
\usepackage{amsfonts}
\usepackage{amssymb}
\usepackage{amsthm}
\usepackage{thmtools}
\usepackage{complexity}
\usepackage[capitalize,nameinlink]{cleveref} 
\crefname{equation}{Equation}{Equations}
\crefname{figure}{Figure}{Figures}
\newtheorem{theorem}{Theorem}[section]
\newtheorem{proposition}[theorem]{Proposition}
\newtheorem{lemma}[theorem]{Lemma}
\newtheorem{corollary}[theorem]{Corollary}

\newtheorem{definition}[theorem]{Definition}

\newtheorem{remark}[theorem]{Remark}
\newtheorem{problem}[theorem]{Problem}
\usepackage{float}
\usepackage{tikz}
\usetikzlibrary{math}
\usetikzlibrary{calc}
\usetikzlibrary{positioning}
\tikzstyle{vtx}=[circle, draw, fill=black, inner sep=0pt, minimum width=5pt]
\tikzstyle{vtx-white}=[circle, draw, fill=white, inner sep=0pt, minimum width=5pt]

\def\ni{\noindent}
\newcommand{\nbc}{neighborhood-balanced $k$-coloring\xspace}

\newcommand{\nbcs}{neighborhood-balanced $k$-colorings\xspace}
\newcommand{\nbcd}{neighborhood-balanced $k$-colored\xspace}

\newcommand{\probNBC}{\textsc{\nbc}\xspace}

\newcommand{\probNBCshort}{\textsc{$k$-NBC}\xspace}

\newcommand{\yes}{\textsf{YES}\xspace}

\newif\ifshownotes
\shownotestrue 

\ifshownotes
  \usepackage[textsize=scriptsize, colorinlistoftodos]{todonotes}
\else
  \usepackage[disable]{todonotes}
\fi

\title{Neighborhood Balanced $k$-Coloring of Graphs}

\authorentry{Maurice Genevieva Almeida}{1}{p20230078@goa.bits-pilani.ac.in}{}{}
\authorentry{Siddharth Gupta}{2}{siddharthg@goa.bits-pilani.ac.in}{0000-0003-4671-9822}{https://guptasid.bitbucket.io/}  
\authorentry{Ravindra Pawar }{3}{cs26r001@cse.iitm.ac.in}{0000-0003-4541-8462}{https://sites.google.com/view/ravi-pawar}
\authorentry{Tarkeshwar Singh}{4}{tksingh@goa.bits-pilani.ac.in}{}{}

\affilentry{1,4}{Department of Mathematics,\\ BITS Pilani, K K Birla Goa Campus, Goa, India.}
\affilentry{2}{Department of Computer Science and Information Systems,\\ BITS Pilani, K K Birla Goa Campus, Goa, India.}
\affilentry{3}{The Institute of Mathematical Sciences, HBNI, Chennai, India.}

\begin{document}

\date{}
\maketitle
\begin{abstract}
For a simple graph $G = (V, E)$ and a positive integer $k \ge 2$, a coloring of vertices of $G$ using exactly $k$ colors such that each vertex has an equal number of neighbors of each color is called {\it neighborhood-balanced $k$-coloring}, and the graph is called \nbcd graph.
This generalizes the notion of neighborhood balanced coloring of graphs introduced by Bryan Freyberg and Alison Marr (Graphs and Combinatorics, 2024).
We derive some necessary/sufficient conditions for a graph to admit a \nbc and discuss several graph classes that admit such colorings.
\end{abstract}

\footnote{Author 1 is the corresponding author of the article}
\section{Introduction}
Let $G=(V,E)$ be a simple graph.
For any vertex $v\in V$, define the \textit{neighborhood} of $v$ as a set $N(v):=\{u: uv\in E\}$.
The members of $N(v)$ are called the \textit{neighbors} of $v$ and the cardinality of $N(v)$ is called the degree of $v$, denoted as $d(v)$.
Further, if we add the vertex $v$ to its neighborhood $N(v)$, then we get the closed neighborhood of $v$, denoted by $N[v]$.
For graph-theoretic notation, we refer to Chartrand and Lesniak \cite{chart}.


Freyberg and Marr \cite{nbc} introduced the concept {\it neighborhood balanced coloring}.
A  neighborhood balanced coloring of a graph $G$ is a vertex coloring of $G$ using two colors, say red and green, such that each vertex has an equal number of neighbors of both colors.
It is easy to see that if a graph admits a neighborhood balanced coloring, then the degree of every vertex is even.
This notion of coloring is somewhat similar to {\it cordial labeling} \cite{cordial}.
In \cite{nbc}, the authors characterized the $2$-regular graphs, complete graphs, and complete multipartite graphs that have a neighborhood balanced coloring.
In addition, they showed that in a neighborhood balanced colored graph $G$, the number of red-green edges is half the number of edges of $G$, and the number of red-red or green-green edges is one-fourth the number of edges of $G$.

In this article, we generalize neighborhood balanced coloring using $k$ colors, called \nbc of graphs, for any positive integer $k \ge 2$.
Formally, a \nbc of a graph is defined as follows.
\begin{definition} \label{def:nbc}
    Let $G$ be a graph, and let $k \ge 2$ be an integer.
    If the vertices of $G$ can be colored using the $k$ colors, say $1, 2,\dots, k$, so that every vertex has an equal number of neighbors of each color, then the coloring is a \emph{\nbc} of $G$.
    A graph that admits such a coloring is called a \emph{\nbcd} graph.
\end{definition}

In other words, a \nbc of graph $G$ is a partition of the vertex set $V$ of $G$ into $k$ sets $V_1, V_2, \dots, V_k$ (that is, color classes) such that every vertex has an equal number of neighbors in each set.
Sometimes it is convenient to use $c\colon V\rightarrow \{-k, \dots, 0, \dots, k\}$ to denote a $(2k+1)$-coloring and $c\colon V\rightarrow \{-k, \dots,-2, -1,1, 2, \dots, k\}$ to denote a $2k$-coloring. 
Let $w(v)=\sum_{u\in N(v)}c(u)$.
Thus, G admits a \nbc if and only if $w(v)=0$ for all $v \in V(G)$.
Note that if $c$ is a \nbc of graph $G$ using colors in the order $(1, 2, \dots, k)$, then the colorings $c=c_1, c_2, \dots, c_k$ obtained by rotating the colors in a cyclic order are also neighborhood-balanced $k$-colorings of $G$.
We refer to such colorings as the colorings obtained by cyclic shift of $c$.

\begin{definition}[Equally colored set] \label{def:equally-colored-set}
Let $c$ be a vertex coloring of a graph $G$.
A subset $X \subseteq V(G)$ is said to be \emph{equally colored} if each color appears the same number of times in $X$.
An equally colored set $X$ is called a rainbow if each color appears in $X$ exactly once.
\end{definition}
Now we can restate the definition of \nbc graph as: a graph $G$ is \nbcd if a set $N(v)$ is equally colored for all vertices $v \in V(G)$.
Refer Figure \ref{fig: nbc of c8} for examples of \nbcd graphs.
As a brief remark, note that the class of \nbcd graphs is not hereditary; this follows readily from the results in \cite{collin-cnbc}.

\bigskip
As an application, a \nbcd graph can describe settings in need of a combination of two or more complementary abilities or perspectives. For instance, imagine a learning environment of a team learning where each student is allocated one of three roles: Analyst (red), Communicator (blue), or Innovator (green).
Draw the students as vertices, and join the vertices by an edge if the respective students tend to work together.
A neighborhood balanced $3$-coloring implies that within each student's immediate social circle of students, each of the three roles is represented in a balanced way.
This organization facilitates overall learning by introducing each student to analytical, communicative, and creative aspects through their day-to-day interactions.\\

In Section \ref{sec:NC}, we prove some necessary conditions for a graph to admit a \nbc and we characterize some graph classes that admit such a coloring.
Further,  we prove that there is no forbidden subgraph characterization for the class of \nbcd graphs.
In Sections \ref{sec:operations} and \ref{sec:examples}, we give various ways to construct \nbcd graphs and study some graph classes that admits \nbc.

\begin{figure}
\centering


\begin{tikzpicture}[scale=1.0]

\node[circle, inner sep = 3pt,, fill=red] (c1) at (-1.4444,1.5556) {};
\node[circle, inner sep = 3pt,, fill=red] (c2) at (-0.8000,0.0000) {};
\node[circle, inner sep = 3pt,, fill=green] (c3) at (-1.4444,-1.5556) {};
\node[circle, inner sep = 3pt,, fill=green] (c4) at (-3.0000,-2.2000) {};
\node[circle, inner sep = 3pt,, fill=red] (c5) at (-4.5556,-1.5556) {};
\node[circle, inner sep = 3pt,, fill=red] (c6) at (-5.2000,0.0000) {};
\node[circle, inner sep = 3pt,, fill=green] (c7) at (-4.5556,1.5556) {};
\node[circle, inner sep = 3pt,, fill=green] (c8) at (-3.0000,2.2000) {};

\draw (c1) -- (c2);
\draw (c2) -- (c3);
\draw (c3) -- (c4);
\draw (c4) -- (c5);
\draw (c5) -- (c6);
\draw (c6) -- (c7);
\draw (c7) -- (c8);
\draw (c8) -- (c1);

\node at (-3,3) {\textbf{Cycle $C_{8}$}};

\node[circle, inner sep = 3pt,, fill=red] (a1) at (2,1.5) {};
\node[circle, inner sep = 3pt,, fill=blue] (a2) at (2,0)   {};
\node[circle, inner sep = 3pt,, fill=green] (a3) at (2,-1.5){};

\node[circle, inner sep = 3pt,, fill=red] (b1) at (5,1.5) {};
\node[circle, inner sep = 3pt,, fill=blue] (b2) at (5,0)   {};
\node[circle, inner sep = 3pt,, fill=green] (b3) at (5,-1.5){};

\draw (a1) -- (b1);
\draw (a1) -- (b2);
\draw (a1) -- (b3);
\draw (a2) -- (b1);
\draw (a2) -- (b2);
\draw (a2) -- (b3);
\draw (a3) -- (b1);
\draw (a3) -- (b2);
\draw (a3) -- (b3);

\node at (3.5,3) {\textbf{Complete bipartite $K_{3,3}$}};


\end{tikzpicture}
\caption{A neighborhood-balanced $2$-coloring of $C_{8}$ and a neighborhood-balanced $3$-coloring of~$K_{3,3}$.}
\label{fig: nbc of c8}
\end{figure}

\section{Necessary Conditions for Neighborhood Balanced $k$-Colored Graphs} \label{sec:NC}
We begin this section by deriving some necessary conditions for graphs to admit a \nbc.
Note that if a graph admits such a coloring, then its minimum degree is $k$.
We state this necessary condition in the following lemma in a more general form.

\begin{lemma}\label{th:prep1}
If a graph $G$ admits a \nbc, then the degree of every vertex is a multiple of $k$. 
\end{lemma}

Observe that a \nbcd graph without isolated vertices cannot have exactly one vertex of any color.

\begin{lemma} \label{th:order of nbc graph}
    If $G$ is a \nbcd graph without isolated vertices, then the order of $G$ is at least $2k$.
\end{lemma}
The lower bound on the order of a \nbcd graph given in the above lemma is sharp.
There are \nbcd graphs of order $2k$, for example, the complete bipartite graph $K_{k,k}$ (refer Theorem \ref{th:multipartite}).

The following lemma shows that if a graph $G$ is \nbcd, where $k$ is even, then it is also neighborhood-balanced $2$-colored.
However, by Lemma \ref{th:prep1}, the converse need not be true.

\begin{lemma}
If a graph $G$ admits a \nbc and $p$ divides $k$, then the graph $G$ also admits a neighborhood-balanced $p$-coloring.
\end{lemma}
\begin{proof}
   Let \( f \colon V(G) \to \{1, 2, \dots, k\} \) be a \nbc\ of \(G\). Define a new coloring by replacing each color \( i \in \{1, 2, \dots, k\} \) with the color \( j \in \{1, 2, \dots, p\} \), where \(j \equiv i\ (\bmod\ p )\). Since \( p \) divides \( k \), each of the \( p \) colors replaces exactly the same number of original colors. Moreover, because every vertex in \(G\) has an equal number of neighbors in each of the original \( k \) color classes, it follows that each vertex still has an equal number of neighbors in each of the \( p \) new color classes.
\end{proof}

For a graph $G$ and subsets $X , Y \subseteq V(G)$, we denote the set of edges joining a vertex in $X$ to a vertex in $Y$ by $E(X,Y)$.
In particular, if $X = Y$, then we write $E(X)$ instead of $E(X,X)$.
Also, we denote a subgraph induced by $X \subseteq V(G)$ by $G\left<X\right>$.
Throughout this article, unless stated otherwise, we assume that $c$ is a $k$-coloring of the graph under consideration, where $k \ge 2$, using colors $1, 2,\dots, k$.
We denote the corresponding color classes as $V_1^{c}(G), V_2^{c}(G),\dots, V_k^{c}(G)$  which form a partition of the vertex set $V(G)$ into $k$ disjoint subsets. If the coloring $c$ and the graph $G$ are clear from the context, then we simply write $V_1, V_2,\dots, V_k$.

\begin{theorem}\label{th:thm1}
    If a graph $G$ admits a \nbc $c$, then
    $$|E(V_i, V_j)|=\frac{2|E(G)|}{k^2}\ \text{and }|E(V_i, V_i)|=\frac{|E(G)|}{k^2}. $$
\end{theorem}
\begin{proof}
    For $i \ne j$, consider a bipartite subgraph $H$ of $G$, induced by the set of edges $E(V_i, V_j)$ with bipartition $V_i \cup V_j$.
    Then for each vertex $v\in V(H)$, $d_H(v) = \frac{1}{k}d_G(v)$.
    The number of edges in $H$ can be counted as
    \begin{equation}\label{eq: edges}
      |E(H)| = |E(V_i, V_j)|=\sum_{v\in V_i}\frac{1}{k}d_G(v)=\sum_{v\in V_j}\frac{1}{k}d_G(v).
    \end{equation}
    This implies $$\sum_{v\in V_i}d_G(v)=\sum_{v\in V_j}d_G(v).$$
    Thus
    \begin{align*}
      2|E(G)|&=\sum_{v\in V_1}d_G(v) + \dots + \sum_{v\in V_{k}}d_G(v)\\
      &=k\sum_{v \in V_i}d_G(v)\\
      &=k^2|E(V_i, V_j)|.
    \end{align*}
    This proves the first equality.
    For the second equality, let $G\left<V_i\right>$ be the subgraph induced by $V_i$.
    For all $v\in V_i$, $d_{G \left< V_i\right>}(v)=\frac{1}{k}d_G(v)$.
    Therefore, by Equation (\ref{eq: edges}), and the first equality of the theorem, we obtain
    \begin{align*}
        2|E(V_i, V_i)|&= \sum_{v \in V_i} d_{G\left<V_i\right>}(v) = \sum_{v \in V_i}\frac{1}{k}d_G(v) =|E(V_i, V_j)|\\
         &=\frac{2|E(G)|}{k^2}
    \end{align*}
    This proves the second equality and completes the proof.
\end{proof}

Note that in general, the color classes in a \nbcd $G$ need not be of the same cardinality.
In Section \ref{sec:unequalcolorclass}, we give a way to construct \nbcd graphs having color classes of different cardinality.
However, in a regular \nbcd graph, all color classes have the same cardinality as stated in the following corollary.

\begin{corollary}\label{th: regular graphs necessary condition}
If $G$ is an $r$-regular graph that admits a \nbc, then $|V_i|=\dfrac{|V(G)|}{k}$.
\end{corollary}
\begin{proof}
For any $i \ne j$, from Equation (\ref{eq: edges}), we have
    \begin{align*}
        |E(V_i, V_j)| &= \frac{1}{k} \sum_{v \in V_i} d_G(v) = \frac{r |V_i|}{k}\\
                      &= \frac{1}{k} \sum_{v \in V_j} d_G(v)=\frac{r |V_j|}{k}.
    \end{align*}
This implies $|V_i| = |V_j|$.
As $|V_1| +|V_2|+  \dots + |V_k| = |V(G)|$, we have $|V_i|=\dfrac{|V(G)|}{k}$ for each $1 \le i \le k$.
\end{proof}
    
\begin{corollary}\label{order and size condition}
    If $G$ is an $r$-regular graph that admits a \nbc, then $|V(G)| \equiv 0\ (\bmod\ k)$ and $|E(G)| \equiv 0\ (\bmod\ {k^2})$.
\end{corollary}

\begin{proof}
    Corollary \ref{th: regular graphs necessary condition} gives $|V(G)|$ is a multiple of $k$ and Theorem \ref{th:thm1} gives $|E(G)|=k^2|E(V_i, V_i)|$.
    This proves the corollary.
\end{proof}

\section{Constructions of Neighborhood-Balanced $k$-Colored Graphs}\label{sec:operations}
We establish necessary conditions under which the {\it product} of two \nbcd graphs is a \nbcd graph.
For completeness, we begin by recalling the definitions of several standard graph products.
Let $G$ and $H$ be graphs.
The {\it cartesian product}  of $G$ and $H$, denoted by $G\ \square \ H$, is a graph with vertex set $V(G)\times V(H)$ and two vertices $(u,v)$ and $(u',v')$ are adjacent if and only if $u=u'$ and $vv'\in E(H)$  or $v=v'$ and $uu'\in E(G)$.
The {\it lexicographic product} of $G$ and $H$, denoted by $G[H]$, is a graph with vertex set $V(G)\times V(H)$ and two vertices $(u,v)$ and $(u',v')$ are adjacent if and only if either $uu'\in E(G)$ or $u=u'$ and $vv'\in E(H)$.
It may be instructive to instead construct $G[H]$ by replacing every vertex of $G$ with a copy of $H$ and then replacing each edge of $G$ with a complete bipartite graph between the corresponding copies of $H$.
The {\it direct product} of $G$ and $H$, denoted by $G\times H$, is a graph with vertex set $V(G)\times V(H)$ and two vertices $(u,v)$ and $(u',v')$ are adjacent if and only if $uu'\in E(G)$ and $vv'\in E(H)$.
The {\it strong product} of $G$ and $H$, denoted by $G \boxtimes H$, is a graph with vertex set $V(G)\times V(H)$ and two vertices $(u,v)$ and $(u',v')$ are adjacent if and only if $u=u'$ and $vv'\in E(H)$ or $v=v'$ and $uu'\in E(G)$ or $uu'\in E(G)$ and $vv'\in E(H)$.
    One can clearly see that the strong product is the edge-disjoint union of the direct product and the cartesian product.
    
For any of the above products and a fixed vertex $u$ of $G$, the set of vertices $\{(u,v):\ v\in V(H)\}$ is called an $H$-layer and we denote it by $H_u$. Similarly, if $v\in V(H)$ is fixed, then the set of vertices $\{(u,v):\ u\in V(G)\}$ is called a $G$-layer and we denote it by $G_v$.
If one constructs $V(G)\times V(H)$ in a natural way, the $H$-layers are represented by rows and the $G$-layers are represented by columns.
We further recall the definition of the join of graphs.
The {\it join of graphs} of $G$ and $H$, denoted by $G+H$, is a graph having vertex set $V(G)\cup V(H)$ and edge set $E(G)\cup E(H)\cup \{xy:  x\in V(G) \text{ and } y \in V(H)\}$.

\begin{theorem}
    If both graphs $G$ and $H$ admit \nbc, then their lexicographic product $G[H]$ admits a \nbc. 
\end{theorem}
\begin{proof}
    Let $c$ and $c'$ be \nbcs of $H$ and $G$, respectively.
    Since the vertex set of $G[H]$ is the union $\cup_{x\in V(G)} V(H_x)$, it is sufficient to specify the coloring scheme for each $H_x$.
    Let the colors used by the \nbc $c$ of $H$ be in the order $(1, 2, \dots, k)$.
    Then, we know that for each $1 \le i \le k$, the coloring $c_i$ obtained by $i$th cyclic shift of $c$ is a \nbc of $H$ that uses colors in the order $(i, i+1, \dots, k, 1, \dots, i-1)$. 
    Since for any $u \in V(G)$, $H_u$ is an isomorphic copy of $H$ in $G[H]$, any \nbc of $H$ can be viewed as a \nbc of $H_u$.
    Therefore, for any vertex $u \in V(G)$ colored $i$ under $c'$, we color $V(H_u)$ using $c_i$.
    Note that since $G$ admits a \nbc $c'$, each vertex $u$ of $G$ is in a unique color class $V_i^{c'}(G)$ and we have colored the $V(H_u)$ by $c_i$.
    Hence, in this way we have colored $V(H_u)$ for all $u \in V(G)$, and thus $V(G[H])$.
    We claim that this is the required coloring scheme.
    
    Let $V_1, V_2, \dots, V_k$ be the partition of $V(G[H])$ into color classes induced by the coloring scheme discussed above.
    Consider a vertex $(u,v) \in V(G[H])$.
    For an arbitrary color $i$, we count the neighbors of $(u,v)$ in $G[H]$ that are colored $i$, that is $|N(u,v) \cap V_i|$.
    Since $H_u$ admits a \nbc, and $(u,v) \in V(H_u)$, $(u,v)$ has an equal number of neighbors in each color class of $H_u$.
    Now let us count the number of neighbors of $(u,v)$ in $V(G[H]) - V(H_u)$ that are colored $i$.
    Since $c'$ is a neighborhood-balanced $k$-coloring of $G$, it follows that vertex $u$ has precisely $p$ neighbors of each color in $G$, for some integer $p$.
    Then, as per our coloring scheme for $G[H]$, for the neighbors of $u$ in the color class $V_j^{c'}(G)$, we have colored the corresponding $p$ copies of $H$ using $c_j$.
    Therefore, the total number of neighbors of $(u,v)$ in $V(G[H]) - V(H_u)$ having color $i$, is given by $p (|V_i^{c_1}(G[H])| + \dots + |V_i^{c_k}(G[H])|)$.
        
    Now since $c_2$ obtained from $c_1$ by a single cyclic shift of colors $(1, 2, \dots, k)$, we obtain $|V_1^{c_1}(G[H])| = |V_2^{c_2}(G[H])|$.
    In general, we have the following equalities:
        \begin{align*}
            &|V_1^{c_1}(G[H])| = |V_2^{c_2}(G[H])| = \dots = |V_k^{c_k}(G[H])|\\
            &|V_2^{c_1}(G[H])| = |V_3^{c_2}(G[H])| = \dots = |V_1^{c_k}(G[H])|\\
            &|V_3^{c_1}(G[H])| = |V_4^{c_2}(G[H])| = \dots = |V_2^{c_k}(G[H])|\\
            &\vdots\\
            &|V_k^{c_1}(G[H])| = |V_1^{c_2}(G[H])| = \dots = |V_{k-1}^{c_k}(G[H])|.
        \end{align*}

    If we substitute $\ell_i = |V_i^{c_1}(G[H])|$, then we obtain that the total number of neighbors of $(u,v)$ in $V(G[H]) - V(H_u)$ that are colored $i$ is $p(\ell_1 + \dots + \ell_k)$ which is a constant.
    Since the color $i$ and vertex $(u,v)$ were arbitrary, we conclude that each vertex of $G[H]$ has an equal number of neighbors in each color class.
    This completes the proof.
\end{proof}

\begin{theorem}
    Let $G$ and $H$ be two graphs. If $H$ admits a \nbc $c$ with $|V^c_i(H)|=|V^c_j(H)|$, then the lexicographic product $G[H]$ admits a \nbc.
\end{theorem}
\begin{proof}
    Apply the \nbc scheme $c$  of the graph $H$ to each $H$-layer in $G[H]$. We claim that this is a \nbc of $G[H]$.\\
    Indeed, let $(u,v)\in V(G[H])$.
    Since $c$ is a \nbc of $H$, we know $(u,v)$ has an equal number of neighbors of each color within the $H$-layer in which it lies.
    Outside of this copy, $(u,v)$ has $d_{G}(u) \times |V^{c}_{i}(H)|$ neighbors with color $i$. Since  $|V^c_i(H)|=|V^c_j(H)|$, the claim follows.
 \end{proof}
 
 \begin{theorem}\label{directproduct}
     If one of the graphs $G$ or $H$ admits a \nbc, then so does the direct product $G\times H$.
 \end{theorem}
 \begin{proof}
     Without loss of generality, assume $G$ admits a \nbc $g$. Consider the graph $G\times H$. Color each $ G$-layer (i.e., column of vertices) using the coloring scheme $g$.
     A vertex $(u,v)\in V(G\times H)$ is adjacent to neighbors of $u$ in those $G$-layers which are due to neighbors of $v$ in $H$. Since all the $G$-layers are \nbcd , $(u,v)$ will have an equal number of neighbors of each of the $k$ colors, and the result follows.
 \end{proof}


\begin{theorem} \label{th:cartesian}
    If both the graphs $G$ and $H$ admit \nbc, then so does the cartesian product $G\, \square\, H$. 
\end{theorem}
\begin{proof}
    Let $g=g_1$ and $h=h_1$ be \nbcs of $G$ and $H$, respectively.
    Let $h_1$ be the coloring that uses colors in the order $(1, 2, \dots, k)$ on the vertices of $H$.
    Then, we know that for each $1 \le i \le k$, $h_i$ is a \nbc of $H$ obtained by $i$th cyclic shift of $h$ and $g_i$ is a \nbc of $G$ obtained by $i$th cyclic shift of $g$.
    
    Consider the graph $G\, \square\, H$.
    Color each $H$-layer (that is, a row of vertices) according to $h=h_1$.
    Then recolor the first $G$-layer (that is, the first column of vertices) according to $g=g_1$, and if the vertex in row $i$ changes color from $1$ to $i$, apply $h_i$ to that row.
    This also ensures that every $G$-layer is colored using one of the colorings $g_1, g_2, \dots, g_{k}$.
    Now we have a coloring of $G\, \square\, H$ in which every $H$-layer has been colored using one of the colorings $h_1, h_2, \dots, h_{k}$ and every $G$-layer has been colored using one of the colorings $g_1, g_2, \dots, g_{k}$.
    As $g_1, g_2, \dots, g_{k}$ and $h_1, h_2, \dots, h_{k}$ are \nbcs, every vertex in $G\, \square\, H$ has an equal number of neighbors of each color in each $G$-layer and $H$-layer.
    This completes the proof.
\end{proof}

The converse of the above theorem is not true. In Theorem \ref{th: hamming graphs}, we show that the hamming graph $H(k,k)$, which is the cartesian product of $k$ complete graphs $K_k$ is  a \nbcd graph but the complete graph $K_k$ is not a \nbcd graph (see Lemma \ref{completegraph}).
\begin{theorem}
    If $G$ and $H$ both admit \nbc, then so does the strong product $G \boxtimes H$. 
\end{theorem}
\begin{proof}
     Let $g=g_1$ and $h=h_1$ be \nbcs of $G$ and $H$, respectively.
    Let $h_1$ be the coloring that uses colors in the order $(1, 2,  \dots, k)$ on the vertices of $H$.
    Then, we know that for each $1 \le i \le k$, $h_i$ is a \nbc of $H$ obtained by $i$th cyclic shift of $h$ and $g_i$ is a \nbc of $G$ obtained by $i$th cyclic shift of $g$.
    Consider the subgraph $G \square H$ of $G \boxtimes H$.
    Color each $H$-layer (that is, a row of vertices) according to $h=h_1$.
    Then recolor the first $G$-layer (that is, the first column of vertices) according to $g=g_1$, and if the vertex in row $i$ changes color from $1$ to $i$, apply $h_i$ to that row.
    In Theorem \ref{th:cartesian}, we proved that if the vertices are colored using such a coloring, then the cartesian product $G \square H$ is \nbcd.
    Now as $g_1, g_2,\dots,g_{k}$ are \nbcs of $G$, each $G$-layer is \nbcd  and hence the direct product $G \times H$ is \nbcd.
   As the strong product is the edge-disjoint union of the direct product and the cartesian product, we have $G \boxtimes H$, neighborhood-balanced $k$-colored.  
\end{proof}

\begin{theorem} \label{th:join}
    If $G$ admits a \nbc $g$ with 
    $|V^g_1(G)|=|V^g_2(G)|=\dots=|V^g_k(G)|$ and $H$ admits a \nbc $h$ with $|V^h_1(H)|=|V^h_2(H)|=\dots=|V^h_k(H)|$, then $G+H$ admits a \nbc.
\end{theorem}
\begin{proof}
  In $G+H$, each vertex of $G$ is adjacent to every vertex of $H$.
  Color the vertices of $G$ using $g$ and color the vertices of $H$  using $h$.
  Consider a vertex $v\in V(G) \subseteq V(G+H)$.
  As $g$ is \nbc of $G$, $v$ has an equal number of neighbors of each color in $G$.
  Also, as $|V^h_1(H)|=|V^h_2(H)|=\dots=|V^h_k(H)|$, $v$ continues to have an equal number of neighbors of each color in $G+H$.
  The same argument works for a vertex in $H$.
  Therefore, $G+H$ is a neighborhood-balanced k-colored graph.
\end{proof}

\ni The following corollary follows from Theorem \ref{th:join} and Corollary \ref{th: regular graphs necessary condition}.

\begin{corollary}
    If $G$ and $H$ are both regular graphs admitting \nbc, then so does $G+H$.
\end{corollary}

We next investigate the \nbcd graphs obtained through unary operation--`union of a graph with itself'.
Note that the disjoint union case is straightforward: the vertex disjoint union of \nbcd graphs is \nbcd.
Therefore, we focus on the the union of a graph with itself, where vertex sets `overlap partially'.
Now we formally define this.
\par
Let $G$ be a graph and let $H$ be a proper subgraph of $G$. 
The {\it union of $n$ copies of $G$ over $H$} is the graph obtained by taking $n$ vertex-disjoint copies of $G$ and identifying the corresponding vertices of their subgraphs isomorphic to $H$.
Similarly, for a nonempty proper subset $S \subsetneq V(G)$, the {\it union of $n$ copies of $G$ over $S$} is defined as the union of $n$ copies of $G$ over the subgraph $G\left< S \right>$ induced by $S$.
We denote this graph as $nG_S$.
A subset $S$ of vertices of a graph is said to be {\it independent} if no two vertices in $S$ are adjacent. 
Otherwise, $S$ is called {\it dependent}.
\par
Our focus is particularly on the study of \nbc of unions over induced subgraphs instead of the union over any arbitrary subgraphs.
If $G$ is a \nbcd and $S \subsetneq V(G)$ is independent, then it is easy to show that $nG_S$ is a \nbcd graph as proved in Theorem \ref{th:nGs}.
However, when $S$ is dependent, the situation is more complicated, as it depends on the integer $n$ as well as on the structure of the subgraph $G\left< S \right>$ induced by $S$.

\begin{theorem}\label{th:nGs}
The union of $n$-copies of a \nbcd graph over an independent set is a \nbcd graph. 
\end{theorem}
\begin{proof}
Let $S \subsetneq V(G)$ be an independent set of vertices in a neighborhood-balanced $k$-colored graph $G$ with corresponding coloring $c$.
Color the vertices of one copy of $G$ in $nG_S$ using the coloring $c$.
For the remaining $(n - 1)$ copies, color all vertices except those in $S$ using the same coloring $c$.
This coloring ensures that all the vertices of $nG_S$, except those of the set $S$, have an equal number of neighbors of each color.
Now, consider any vertex $v \in S$.
Since $S$ is an independent set in $G$, it follows that $\deg_{nG_S}(v) = n \cdot \deg_G(v)$.
Moreover, as $c$ is a neighborhood-balanced $k$-coloring of $G$, $v$ has $\frac{\deg_G(v)}{k}$ neighbors of every color in $G$.
Therefore, in $nG_S$, the vertex $v$ has $\frac{n \cdot \deg_G(v)}{k}$ neighbors of each color.
\end{proof}

Next, we first derive a necessary condition on the integer $n$ for a graph $nG_S$ to admit a \nbc when $S \subsetneq V(G)$ is dependent.

\begin{theorem}\label{th: unionoversepset}
If the union of $n$-copies of \nbcd graph over a dependent set is \nbcd then $n\equiv 1 \bmod(\mathrm{lcm}\ (L,M))$, where $L = \mathrm{lcm}\ \left( \frac{k}{\gcd(q_1, k)}, \frac{k}{\gcd(q_2, k)}, \dots, \frac{k}{\gcd(q_r, k)} \right)$, $M = \frac{k^2}{\gcd(p, k^2)}$, where $q_i = \deg_{G\langle S \rangle}(v_i)$,  where $1 \le i \le r$, and  $p = |E(G\langle S \rangle)|$.
\end{theorem}
\begin{proof}
    Note that $\deg_{nG_S}(v_i) = n(\deg_G(v_i)-q_i) + q_i$, for $1 \leq i \leq r$.
    Since both $G$ and $nG_S$ are \nbcd graphs, by Lemma \ref{th:prep1}, we get $nq_i \equiv q_i\pmod{k},$ for $1\leq i \leq r$.
    Again since $nG_S$ is a \nbcd graph, by Theorem \ref{th:thm1}, $|E(nG_S)| \equiv 0 \pmod{k^2}$.
    Now $|E(nG_S)|= n(|E(G)|-p) + p$ and as $G$ is \nbcd graph, we get $np \equiv p \pmod{k^2}$. 
    This leads us to the following system of linear congruences:
    \begin{align*}
        nq_i &\equiv q_i\ (\bmod\ k) \quad (1\leq i \leq r), \\
        np &\equiv p\ (\bmod\ k^2).
    \end{align*}
    That is, 
    \begin{align*}
        (n-1)q_i &\equiv 0\ (\bmod\ k)\quad  (1 \leq i \leq r), \\
        (n-1)p &\equiv 0\ (\bmod\ k^2).
    \end{align*}
    This implies, 
     \begin{align*}
        (n-1) &\equiv 0\ (\bmod\ \tfrac{k}{\gcd(k,q_i)}) \quad  (1 \leq i \leq r) ,\\
        (n-1) &\equiv 0\ (\bmod\ \tfrac{k^2}{\gcd(k^2,p)}).
    \end{align*}
    The first $r$-congruence equations in the above system imply
    \begin{align*}
        (n-1)&\equiv 0\ (\bmod\ L),\ \text{where } L = \mathrm{lcm} \big(\tfrac{k}{\gcd(k, q_1)}, \tfrac{k}{\gcd(k, q_2)},  \dots,\tfrac{k}{ \gcd(k, q_r)}\big),\\
            (n-1)& \equiv 0\ (\bmod\ M),\ \text{where } M = \tfrac{k^2}{\gcd(k^2, p)}.
    \end{align*}
That is, \begin{align*}
    n&\equiv 1\ (\bmod\ L),\\
    n&\equiv 1\ (\bmod\ M).
\end{align*}
This system has a solution if and only if $\mathrm{gcd}(L,M)$ divides $(1-1)=0$, which is always true and the solution is $n\equiv 1\ (\bmod\ \mathrm {lcm}(L,M)).$
\end{proof}

\begin{corollary}
If the one edge union of $n$-copies of a \nbcd graph is \nbcd, then $n\equiv 1 (\bmod\ k^2)$. 
\end{corollary}

Next, we show that the necessary condition on the integer $n$ obtained in Theorem \ref{th: unionoversepset} may or may not be sufficient.
For some graphs, the converse may depend on the structure of the subgraph induced by the dependent set under consideration, as seen in the following example.\par
Let the terminology be as used in Theorem \ref{th: unionoversepset}.
For $m\equiv 0\ (\bmod\ 4)$, Freyberg and Marr \cite{nbc}, showed that $C_m$ is a neighborhood-balanced $2$-colored graph.
We now study the union of $n$-copies of $C_m$ for $m \equiv 0\ (\bmod\ 4)$, over a nonempty proper dependent set $S$.
Note that $C_m\left< S \right>$ is a disjoint union of paths.
If $C_m\left< S \right>$ contains an odd number of edges, then $p$ is odd and $M=\frac{4}{\gcd(p,4)}=4$.
Further, as $S$ is a proper subset of $V(C_m)$, at least two $q_i$ must be equal to $1$.
This implies at least two of $\gcd(q_i,2)$ will be equal to $1$ and hence at least two of $\frac{2}{\gcd(q_i,2)}=2$.
This implies $L=\mathrm{lcm}\big(\frac{2}{\gcd(q_1,2)}, \frac{2}{\gcd(q_2,2)}, \dots, \frac{2}{\gcd(q_r,2)}\big)=2$.
Thus $\mathrm{lcm}(L,M) = \mathrm{lcm}(2,4)=4$ and so $n \equiv 1 \pmod{4}$.
If $C_m\left< S \right>$ contains an even number of edges, then $p$ is even. 
Therefore, $M=\frac{4}{\gcd(p,4)}=1\ \text{or}\ 2$.
As argued in the above paragraph, we have $L = 2$.
Thus $\mathrm{lcm}(L,M)=2$ and so $n\equiv 1 (\bmod\ 2)$.\\

Next, we introduce the concept of an ideal dependent set in $V(C_m)$.
A nonempty proper dependent set $S \subsetneq V(C_m)$ is said to be {\it ideal dependent set} if the subgraph $C_m\left< S \right>$ induced by $S$ has the following properties:
\begin{itemize}
\setlength\itemsep{0pt}
    \item[i.] It has no trivial (consisting of exactly one vertex) components.
    \item[ii.] The number of vertices in $C_m$ that appear between any two consecutive components is odd.
    \item[iii.] If $C_m\left< S \right>$ has only one component, then it is a path of even length (such a path is called an even path).
\end{itemize}
\begin{remark}\label{remark on components}
    If $S$ is an ideal dependent set in $C_m$ ($m$ even) such that $C_m\left<S\right>$ has only one component (i.e. an even path), then there is an odd number of vertices between the end vertices of the path. So, if we consider this lone component of $C_m\left<S\right>$ as two consecutive (repeated) components along the cycle $C_m$, then there is an odd number of vertices between them. Similarly, if $S$ is not an ideal dependent set in $C_m$ ($m$ even) such that $C_m\left<S\right>$ has only one component (i.e. an odd path), then considering the odd path as two consecutive (repeated) components along the cycle $C_m$, there is an even number of vertices between them.
\end{remark}
Figure \ref{fig:segments} denotes the ideal dependent set $S= \{0,1,2,6,7,8\}$ in $C_{10}$.
\begin{figure}
    \centering
    \begin{tikzpicture}
  \foreach \i in {0,...,9} {
    \node[circle, draw, fill=white, inner sep=1.5pt] (v\i) at ({360/10*\i}:2cm) {$\i$};
  }
\draw[thick,blue] (v0) -- (v1);
\draw[thick,blue] (v1) -- (v2);
\draw[] (v2) -- (v3);
\draw[] (v3) -- (v4);
\draw[] (v4) -- (v5);
\draw[] (v5) -- (v6);
\draw[thick,blue] (v6) -- (v7);
\draw[thick,blue] (v7) -- (v8);
\draw[] (v8) -- (v9);
\draw[] (v9) -- (v0);
\end{tikzpicture}
    \caption{Ideal dependent set $S = \{0,1,2,6,7,8\}$ in $C_{10}$ with the two components colored blue.}
    \label{fig:segments}
\end{figure}
We have the following observation about the graph induced by an ideal dependent subset of $V(C_m)$.

\begin{lemma}
    If $S$ is an ideal dependent set in $C_m$ for $m\equiv 0 (\bmod\ 4)$, then $|E(C_m\left< S \right>)|$ is even.
\end{lemma}
\begin{proof}
    If $S$ is such that $C_m\left< S \right>$ has only one component, then $C_m\left< S \right>$ has to be an even path, and the result follows.
    If $C_m\left< S \right>$ has more than one component, then there is an odd number of vertices, and hence an even number of edges between any two consecutive components.
    Therefore, as $C_m$ also has an even number of edges, $|E(C_m\left< S \right>)|$ is even.
\end{proof}
Now, we give a complete characterization of the dependent sets $S \subsetneq V(C_m)$ such that the union of $n$ copies of a cycle $C_m$ for $m \equiv 0 (\bmod\ 4)$ is \nbcd.

\begin{theorem}\label{th: unionoverdepsetcycle}
    The union of $n$-copies of $C_m$ over a dependent set $S$ is a neighborhood-balanced $2$-colored graph if and only if $S$ is an ideal dependent set in $C_m$, where $n$ is odd and  $m\equiv 0 (\bmod\ 4)$.
\end{theorem}
\begin{proof}
    Let $S$ be an ideal dependent set in $C_m$ that induces $\ell$ components and let $G$ be the union of $n$-copies of $C_m$ over $S$.
    Let $V(C_m) = \{v_1, v_2,  \dots, v_m\}$ and for $1 \le i \le n$, let $v_j^i$ be the copy of $v_j$ in the $i$th copy of $C_m$.
    Then the vertex set of $G$ is: $V(G) = S \cup \big(\cup_{i=1}^n \{u^i :  u \notin S\}\big)$.
    Now we partition the set $V(C_m)-S$ as follows.
    For each $1 \le i \le \ell$, let $T_i = \{v_{i,\alpha_{i}+1}, v_{i,\alpha_{i}+2}, \dots, v_{i,\alpha_{i} + t_i}\}$ be the set of consecutive vertices of $C_m$ that appears between two consecutive components $C_i$ and $C_{i+1}$ (if $C_m\left<S\right>$ consists of only one component then we take $C_i=C_{i+1}$), where the second index in the suffix of a vertex $v$ indicates its position in the cyclic ordering of $V(C_m)$.
    Note that each $t_i$ is odd (refer definition of ideal dependent set and Remark \ref{remark on components}).
    This gives us a partition of $V(G)-S$ as: $T_{i}^{j} = \{v_{i,\alpha_{i}+1}^{j}, v_{i,\alpha_{i}+2}^{j},  \dots, v_{i,\alpha_{i} + t_i}^{j}\}$, where $ 1 \le j \le n$ and $1\leq i \leq \ell$.
   \par
    Let $c$ be a neighborhood-balanced $2$-coloring of $C_m$ and $\overline{c}$ be its cyclic shift.
    We now provide a coloring $c'$ for the graph $G$ as follows.
    \begin{align*}
        c'(v)=
        \begin{cases}
            c(v) &\text{if } v \in S,\\
            c(v_{i, \alpha_{i}+s}) &\text{if } v = v_{i, \alpha_{i} + s}^{j} \text{ for $1 \leq j \leq \frac{n+1}{2}$ and $s\equiv 1 \text{ or }3\  (\bmod\ 4)$},\\
            \overline{c}(v_{i,\alpha_i+s}) &\text{if } v = v_{i, \alpha_{i} + s}^{j}\  \text{for $\frac{n+3}{2}\leq j \leq n$ and $s\equiv 1\text{ or }3\ (\bmod\ 4)$},\\
            c(v_{i,\alpha_i+s}) &\text{if } v=v_{i, \alpha_{i} + s}^{j}\ \text{for $1\leq j \leq n$ and $s\equiv 0\text{ or }2\ (\bmod\ 4)$}.
        \end{cases}
    \end{align*}
    Note that $c'(V(C_m)) = c(V(C_m))$.
    The vertex $v_{\alpha_i}$ had an equal number of neighbors of both colors in the original copy of $C_m$.
    It receives an additional $\frac{n-1}{2}$ neighbors of each color.
    Therefore, it has an equal number of neighbors of both colors in $G$.
    Similarly, as each $t_i$ is odd, using a similar argument, we can show that each vertex $u_{\alpha_i+t_i+1}$ has an equal number of neighbors of both colors in $G$.
    As $c$ and $\overline{c}$ are neighborhood-balanced $2$-colorings, it is easy to check that the vertices $v_{i,\alpha_{i}+1}, v_{i,\alpha_{i}+2}, \dots, v_{i,\alpha_{i} + t_i}$ have one neighbor of each color.
    \par
    Conversely, suppose that $S$ is not an ideal dependent set in $C_m$.
    So the subgraph $C_m\left< S \right>$ is either an odd path or there exist two components such that there is an even number of vertices between them.  
    Let us assume the components to be $C_i$ and $C_{i+1}$ (if $C_m\left< S \right>$ consists of a single component then we take $C_i=C_{i+1}$).
    Let the vertices between these two components be $v^j_{i,\alpha_i+1}, v^j_{i,\alpha_i+2}, \dots, v^j_{i,\alpha_i+t_i}$, where $1\leq j \leq n$ and $t_i$ is even (refer definition of ideal dependent set and Remark \ref{remark on components}).
    Consider the vertex $v_{\alpha_i}$.
    Its $n+1$ neighbors are $\{v_{\alpha_i-1}, v^j_{i,\alpha_i+1} : 1 \leq j \leq n \}$.
    So $\frac{n+1}{2}$ of these must be of one color and the other $\frac{n+1}{2}$ must be of the other color.
    Now consider the vertices $v^j_{i,\alpha_i+2}$ for $1\leq j \leq n$.
    As all the sets $N(v^j_{i,\alpha_i+1})$ for $1\leq j \leq n$ are equally colored and for any fixed $1 \leq p \leq n$, the neighbors of $v^p_{i,\alpha_i+1}$ are $v_{\alpha_i}$ and $v^p_{i,\alpha_i+2}$,  all the vertices of type $v^j_{i,\alpha_i+2}$ receive the same color that is complementary to the color received by $v_{\alpha_i}$.
    Further, for each  $1 \le j \le n$, the vertex $v^j_{i,\alpha_i+4}$ receives the color complementary to that of the $v^j_{i,\alpha_i+2}$.
    Therefore, all the vertices of the type $v^j_{i,\alpha_i+4}$, where $1 \le j \le n$ receive the same color as complementary to the color of the vertices of the type $v^j_{i,\alpha_i+2}$ where $1\leq j \leq n$.
    In a similar way, all the vertices of type $v^j_{i,\alpha_i+6}, v^j_{i,\alpha_i+8}, \dots, v^j_{i,\alpha_i+t_i}$ for $1\leq j \leq n$ receive the same color.
    This, however, implies that the set $N(v_{\alpha_i+t_i+1})$ is not equally colored, as possibly only one of its neighbors can have a different color. 
\end{proof}
Theorem \ref{th: unionoverdepsetcycle} characterizes the dependent set $S$ for which the union of $n$-copies of $C_m$ for $m\equiv 0 \pmod{4}$ over $S$ is \nbcd graph. If for a particular graph there exists no such set $S$, then we can take $S=\phi$. This motivates the following problem.
\begin{problem}
    Given a \nbcd graph $G$ and an integer $n$, satisfying conditions as in Theorem \ref{th: unionoversepset}, characterize the dependent set $S$ (that is, the subgraph $G\left< S \right>$ induced by $S$) such that the union of $G$ over $S$ is a \nbcd graph.
\end{problem}

\begin{theorem} Let $G_1$ and $G_2$ be graphs admitting neighborhood-balanced $k_1$-coloring and neighborhood-balanced $k_2$-coloring, respectively. Then the direct product $G_1\times G_2$ admits a neighborhood-balanced $k_1k_2$-coloring.
\end{theorem}

\begin{proof}
Let
$f:V(G_1)\rightarrow \{1,2,\ldots,k_1\}$
be a neighborhood-balanced $k_1$-coloring of $G_1$, and let
$g:V(G_2)\rightarrow \{1,2,\ldots,k_2\}$
be a neighborhood-balanced $k_2$-coloring of $G_2$.
For each vertex $v\in V(G_2)$, consider the subgraph
\[
G_v=\{(u,v):u\in V(G_1)\},
\]
called the \emph{$G_1$-layer} corresponding to $v$.

\ni Suppose $g(v)=i$, where $1\le i\le k_2$. We color the vertices of the layer $G_v$ according to the coloring $f$ by assigning
\[
h(u,v)=k_1(i-1)+f(u).
\]
Thus, each $G_1$-layer corresponding to a vertex of color $i$ in $G_2$ is colored using the block of colors
\[
\{k_1(i-1)+1,\;k_1(i-1)+2,\;\ldots,\;k_1i\},
\]
following exactly the color pattern of the neighborhood-balanced coloring $f$ of $G_1$.
 Since the colors assigned to different layers of $G_1$ belong to disjoint blocks, the coloring $h$ uses exactly $k_1k_2$ distinct colors.\\
Now let $(x,y)\in V(G_1\times G_2)$, where $g(y)=i$. Since
\[
N_{G_1\times G_2}(x,y)
=
\{(u,v):u\in N_{G_1}(x),\,v\in N_{G_2}(y)\},
\]
every neighbor of $(x,y)$ belongs to a layer corresponding to a neighbor of $y$ in $G_2$.\\
\ni 
Fix an arbitrary color
$
c=k_1(j-1)+r$,
where $1\le j\le k_2$ and $1\le r\le k_1$.
A neighbor of $(x,y)$ has color $c$ if and only if its second coordinate is a neighbor of $y$ colored $j$ in $G_2$, and within the corresponding $G_1$-layer, its first coordinate is a neighbor of $x$ colored $r$ under $f$.

\ni
Since $g$ is a neighborhood-balanced $k_2$-coloring,
\[
\left|N_{G_2}(y)\cap g^{-1}(j)\right|
=
\frac{d_{G_2}(y)}{k_2}.
\]
Similarly, since $f$ is a neighborhood-balanced $k_1$-coloring,
\[
\left|N_{G_1}(x)\cap f^{-1}(r)\right|
=
\frac{d_{G_1}(x)}{k_1}.
\]
\ni
Therefore, the number of neighbors of $(x,y)$ having color $c$ is
\[
\frac{d_{G_2}(y)}{k_2}
\cdot
\frac{d_{G_1}(x)}{k_1}
=
\frac{d_{G_1}(x)d_{G_2}(y)}{k_1k_2}.
\]
\ni
Since
$d_{G_1\times G_2}(x,y)
=
d_{G_1}(x)\,d_{G_2}(y)$,
it follows that every color appears exactly
$\frac{d_{G_1\times G_2}(x,y)}{k_1k_2}$
times in the neighborhood of $(x,y)$.

\ni Hence, every vertex of $G_1\times G_2$ has an equal number of neighbors in each of the $k_1k_2$ colors. Therefore, $h$ is a neighborhood-balanced $k_1k_2$-coloring of $G_1\times G_2$.
\end{proof}

\section{Examples of Neighborhood-Balanced $k$-Colored Graphs} \label{sec:examples}
The converse to Theorem \ref{th:thm1} and hence of Corollary \ref{order and size condition} need not be true in general. For example, the graph $K_{8n}$ satisfies all the necessary conditions mentioned in the Corollary \ref{order and size condition}, but it does not admit a neighborhood-balanced $2$-coloring (see Lemma \ref{completegraph}). 
It would be interesting to study for which graph classes the converse is true.
For example, the converse of Corollary \ref{order and size condition} holds for hamming graphs $H(nk,k)$ (see Theorem \ref{th: hamming graphs}), complete multipartite graphs $K_{n_1,n_2,\dots,n_r}$, $n_i\equiv 0\ (\bmod\ k)$ (see Theorem \ref{th:multipartite}).
In this subsection, we study some regular graph classes that admit a \nbc.
\par
We know that the complete graph $K_n$ is $(n-1)$-regular.
If $K_n$ is a \nbcd graph, then by Lemma \ref{th:prep1}, $n-1$ is a multiple of $k$ and by Corollary \ref{th: regular graphs necessary condition}, $n$ is a multiple of $k$, which is not possible.
We state this observation as a lemma.

\begin{lemma}\label{completegraph}
    A complete graph $K_n$ for $n \geq 2$ does not admit a \nbc.
\end{lemma}

\begin{theorem}\label{th:multipartite}
Let $p \geq 2$.
The complete multipartite graph $K_{n_1,n_2,\dots,n_p}$ admits a \nbc if and only if $n_i\equiv 0\ (\bmod\ k)$ for $i=1,2, \dots, p$.
\end{theorem}
\begin{proof}
    Let $G \cong K_{n_1,n_2,\dots,n_p}$ be a complete multipartite graph. 
    If $n_i \equiv 0\ (\bmod\ k)$, then we can color $\tfrac{n_i}{k}$ vertices using color $i$ for $i = 1, 2, \dots, k$.
    It is easy to see that this is a \nbc of $G$.
    On the other hand, suppose $G$ admits a \nbc. Let $A_1,A_2,\dots, A_p$ be the $p$-partite sets such that $|A_i|=n_i$ and let $A_{i}^{j}$ is the set of vertices in $A_i$ colored $j$.
    Recall that $V_i$ denotes the vertices of $G$ colored $i$.
    Note that for a vertex $v \in V(G)$, there exists some $1 \le \ell \le p$ such that $v \in A_{\ell}$.
    The number of neighbors of such a vertex $v$ colored $i$ is $|V_i| - |A_{\ell}^{i}|$.
    Therefore, the following equation
    \begin{equation}\label{eq: eq2}
        |V_i| - |A_{\ell}^{i}| = |V_j| - |A_{\ell}^{j}|
    \end{equation}
    holds for every $\ell = 1, 2,  \dots, p$.
    Adding these $p$ equations yields $(p-1) |V_i| =(p-1) |V_j|$.
    Thus $|V_i| = |V_j|$.
    It follows from Equation \ref{eq: eq2} that $|A_{\ell}^{i}| = |A_{\ell}^{j}|$ for $\ell = 1, 2, \dots, p$.
    Hence, we conclude that each $n_i$ is a multiple of $k$, and this completes the proof.
\end{proof}

Next, we derive a sufficient condition for some subclasses of circulant graphs to admit a \nbc.
For integers $0 < a_1 < a_2 < \dots <a_k< \frac{n}{2}$, the circulant graph $C_n(a_1,a_2,\dots,a_k)$ is a $2k$-regular graph having vertex set $\mathbb{Z}_n$ with $N(i)=\{i-a_k, i-a_{k-1},\dots,i-a_1, i+a_1, \dots, i+a_{k-1},i+a_k \}$ where all arithmetic is done in $\mathbb{Z}_n$.
\par
Freyberg and Marr \cite{nbc} completely characterized quartic circulant graphs $C_n(a,b)$ that admit a neighborhood-balanced $2$-coloring.

\begin{figure}[htbp]
\centering
\begin{tikzpicture}[scale=2.5, every node/.style={circle, draw, fill=blue!20, minimum size = 5mm, inner sep=1pt}]
  \node (0) at (90:1.5) {0};
  \node (1) at (110:1.5) {1};
  \node (2) at (130:1.5) {-1};
  \node (3) at (150:1.5) {0};
  \node (4) at (170:1.5) {1};
  \node (5) at (190:1.5) {-1};
  \node (6) at (210:1.5) {0};
  \node (7) at (230:1.5) {1};
  \node (8) at (250:1.5) {-1};
  \node (9) at (270:1.5) {0};
  \node (10) at (290:1.5) {1};
  \node (11) at (310:1.5) {-1};
  \node (12) at (330:1.5) {0};
  \node (13) at (350:1.5) {1};
  \node (14) at (10:1.5) {-1};
  \node (15) at (30:1.5) {0};
  \node (16) at (50:1.5) {1};
  \node (17) at (70:1.5) {-1};

  \draw (0) -- (1);
  \draw (1) -- (2);
  \draw (2) -- (3);
  \draw (3) -- (4);
  \draw (4) -- (5);
  \draw (5) -- (6);
  \draw (6) -- (7);
  \draw (7) -- (8);
  \draw (8) -- (9);
  \draw (9) -- (10);
  \draw (10) -- (11);
  \draw (11) -- (12);
  \draw (12) -- (13);
  \draw (13) -- (14);
  \draw (14) -- (15);
  \draw (15) -- (16);
  \draw (16) -- (17);
  \draw (17) -- (0);

  \draw (0) -- (3);
  \draw (1) -- (4);
  \draw (2) -- (5);
  \draw (3) -- (6);
  \draw (4) -- (7);
  \draw (5) -- (8);
  \draw (6) -- (9);
  \draw (7) -- (10);
  \draw (8) -- (11);
  \draw (9) -- (12);
  \draw (10) -- (13);
  \draw (11) -- (14);
  \draw (12) -- (15);
  \draw (13) -- (16);
  \draw (14) -- (17);
  \draw (15) -- (0);
  \draw (16) -- (1);
  \draw (17) -- (2);

  \draw (0) -- (5);
  \draw (1) -- (6);
  \draw (2) -- (7);
  \draw (3) -- (8);
  \draw (4) -- (9);
  \draw (5) -- (10);
  \draw (6) -- (11);
  \draw (7) -- (12);
  \draw (8) -- (13);
  \draw (9) -- (14);
  \draw (10) -- (15);
  \draw (11) -- (16);
  \draw (12) -- (17);
  \draw (13) -- (0);
  \draw (14) -- (1);
  \draw (15) -- (2);
  \draw (16) -- (3);
  \draw (17) -- (4);

\end{tikzpicture}
\caption{A neighborhood-balanced $3$-coloring of $C_{18}(1,3,5)$.}
\label{fig: 3nbc_of_c18(1,3,5)}
\end{figure}

\begin{lemma}\label{circ4k+2}
    Let $a_1, a_2, \dots, a_{2k+1}$ be positive integers such that $1\leq a_1<a_2<\dots < a_{2k+1}<\frac{n}{2}$ and $a_{i+1}-a_{i} \equiv p\ (\bmod\ {2k+1})$, where $n\equiv 0\ (\bmod\ {2k+1})$ and $p \not \equiv 0\ (\bmod\ {2k+1})$.
    Then $C_n(a_1, a_2, \dots, a_{2k+1})$ is neighborhood-balanced $(2k+1)$-colored.
\end{lemma}
\begin{proof}
    Let $G \cong C_n(a_1, a_2, \dots, a_{2k+1})$ with $V(G) = \{0, 1,  \dots, n-1\}$.
    Define a coloring~$c\colon V(G)\rightarrow\{-k, \dots, k\}$ by 
    \begin{align*}
        c(i)=
        \begin{cases}
           \  0 &\text{if $i\equiv 0\ (\bmod\ {2k+1})$,}\\
            \ j &\text{if $i\equiv 2j\ (\bmod\ {2k+1})$; $j=1,2,\dots,k,$}\\
             -j &\text{if $i\equiv 2j-1\ (\bmod\ {2k+1})$; $j=1,2,\dots,k.$}\\
        \end{cases}
    \end{align*}
    Refer Figure \ref{fig: 3nbc_of_c18(1,3,5)} for one such coloring.
    Let $u\in V(G)$ be given. We have, 
    \begin{equation*}
        w(u)=c(u-a_{1})+c(u-a_{2})+\dots+c(u-a_{2k+1})+c(u+a_{1})+c(u+a_{2})+\dots+c(u+a_{2k+1}).
    \end{equation*}
    Suppose that $u+a_1\equiv q\ (\bmod\ {2k+1})$.
    As $a_{i+1}-a_{i}\equiv p\ (\bmod\ {2k+1})$; $p\in \{1,2,\dots,2k\}$, we have $u+a_i\equiv q+(i-1)p\ (\bmod\ {2k+1})$ for $2\leq i \leq 2k+1$.
    As $p\in \{1, 2,  \dots, 2k\}$, $u+a_i$ is congruent to $1, 2, \dots, (2k+1)$ modulo $(2k+1)$ as $i$ takes values from $1, 2,  \dots, (2k+1)$. 
     So $c(u+a_{1})+c(u+a_{2})+\dots+c(u+a_{2k+1})=0$.
     Similar calculations show that $c(u-a_{1})+c(u-a_{2})+\dots+c(u-a_{2k+1})=0$.
     Thus $w(u)=0$ and the coloring $c$ is a neighborhood-balanced $(2k+1)$-coloring of G.
\end{proof}

\begin{figure}[htbp]
\centering
\begin{tikzpicture}[scale=4, every node/.style={circle, draw, fill=blue!20, minimum size = 5mm, inner sep=1pt, text centered}]

\node (0)  at (90:1) {-1};
\node (1)  at (105:1) {-1};
\node (2)  at (120:1) {-1};
\node (3)  at (135:1) {1};
\node (4)  at (150:1) {1};
\node (5)  at (165:1) {1};
\node (6)  at (180:1) {-2};
\node (7)  at (195:1) {-2};
\node (8)  at (210:1) {-2};
\node (9)  at (225:1) {2};
\node (10) at (240:1) {2};
\node (11) at (255:1) {2};
\node (12) at (270:1) {-1};
\node (13) at (285:1) {-1};
\node (14) at (300:1) {-1};
\node (15) at (315:1) {1};
\node (16) at (330:1) {1};
\node (17) at (345:1) {1};
\node (18) at (0:1)   {-2};
\node (19) at (15:1)  {-2};
\node (20) at (30:1)  {-2};
\node (21) at (45:1)  {2};
\node (22) at (60:1)  {2};
\node (23) at (75:1)  {2};

\draw (0) -- (1);   \draw (1) -- (2);   \draw (2) -- (3);
\draw (3) -- (4);   \draw (4) -- (5);   \draw (5) -- (6);
\draw (6) -- (7);   \draw (7) -- (8);   \draw (8) -- (9);
\draw (9) -- (10);  \draw (10) -- (11);\draw (11) -- (12);
\draw (12) -- (13); \draw (13) -- (14);\draw (14) -- (15);
\draw (15) -- (16); \draw (16) -- (17);\draw (17) -- (18);
\draw (18) -- (19); \draw (19) -- (20);\draw (20) -- (21);
\draw (21) -- (22); \draw (22) -- (23);\draw (23) -- (0);

\draw (0) -- (4);   \draw (1) -- (5);   \draw (2) -- (6);
\draw (3) -- (7);   \draw (4) -- (8);   \draw (5) -- (9);
\draw (6) -- (10);  \draw (7) -- (11);  \draw (8) -- (12);
\draw (9) -- (13);  \draw (10) -- (14);\draw (11) -- (15);
\draw (12) -- (16); \draw (13) -- (17);\draw (14) -- (18);
\draw (15) -- (19); \draw (16) -- (20);\draw (17) -- (21);
\draw (18) -- (22); \draw (19) -- (23);\draw (20) -- (0);
\draw (21) -- (1);  \draw (22) -- (2); \draw (23) -- (3);

\draw (0) -- (7);   \draw (1) -- (8);   \draw (2) -- (9);
\draw (3) -- (10);  \draw (4) -- (11);  \draw (5) -- (12);
\draw (6) -- (13);  \draw (7) -- (14);  \draw (8) -- (15);
\draw (9) -- (16);  \draw (10) -- (17);\draw (11) -- (18);
\draw (12) -- (19); \draw (13) -- (20);\draw (14) -- (21);
\draw (15) -- (22); \draw (16) -- (23);\draw (17) -- (0);
\draw (18) -- (1);  \draw (19) -- (2); \draw (20) -- (3);
\draw (21) -- (4);  \draw (22) -- (5); \draw (23) -- (6);

\end{tikzpicture}
\caption{A neighborhood-balanced $4$-coloring of $C_{24}(1,4,7,10)$.}
\label{fig:4nbc_c24}
\end{figure}

\begin{lemma}\label{circ4k}
    Let $a_1,a_2,\dots,a_{2k}$ be positive integers such that $1\leq a_1< a_2 < \dots < a_{2k}<\frac{n}{2}$ and $a_{i+1}-a_{i} \equiv p\ (\bmod\ {2k})$, where $n\equiv 0\ (\bmod\ {2k})$ and $p \not \equiv  0\ (\bmod\ {2k})$.
    Then $C_n(a_1,a_2,\dots,a_{2k})$ is neighborhood-balanced $2k$-colored.
\end{lemma}
\begin{proof}
    Let $G \cong C_n(a_1, a_2, \dots, a_{2k})$ with $V(G)=\{0,1,2,\dots,n-1\}$.
    Define a coloring $c\colon V(G)\rightarrow\{-k, \dots, -1, 1, \dots, k\}$ by 
    \begin{align*}
        c(i)=
        \begin{cases}
           \frac{j}{2}+1 &\text{if } i\equiv jp+1, jp+2, \dots, (j+1)p\ (\bmod\ {2k}),\; j \text{ even,}\\
           - \big(\frac{j+1}{2}\big) &\text{if } i\equiv jp+1, jp+2, \dots, (j+1)p\ (\bmod\ {2k}),\; j \text{ odd.}
        \end{cases}
    \end{align*}
    Refer Figure \ref{fig:4nbc_c24} for one such coloring.
    Let $u\in V(G)$ be given.
    We have, 
    \begin{equation*}
        w(u)=c(u-a_{1}) + c(u-a_{2}) + \dots + c(u-a_{2k}) + c(u+a_{1}) + c(u+a_{2}) + \dots + c(u+a_{2k}).
    \end{equation*}
     As $a_{i+1}-a_{i} \equiv p (\bmod\ {2k})$; $p\in \{1, 2, \dots, 2k-1\}$, under the coloring $c$, the vertices $u+a_1,u+a_2,\dots,u+a_{2k}$ receives each of the colors from the set $\{-k, \dots,-1,1, \dots, k\}$ exactly once.
     So $c(u+a_{1})+c(u+a_{2})+\dots+c(u+a_{2k})=0$.
     Similar calculations show that $c(u-a_{1})+c(u-a_{2})+\dots+c(u-a_{2k})=0$.
     Thus $w(u)=0$ and the coloring $c$ is a neighborhood-balanced $2k$-coloring of G.
\end{proof}

\noindent The proof of the following theorem follows from Lemmas \ref{circ4k+2} and \ref{circ4k}.
\begin{theorem}
    Let $a_1,a_2,\dots,a_k$ be positive integers such that $1\leq a_1 < a_2 < \dots < a_k < \frac{n}{2}$ and $a_{i+1}-a_i \equiv p (\bmod\ {2k})$, where $n\equiv 0 (\bmod\ {2k})$ and $\  p\not\equiv 0 (\bmod\ {2k})$.
    Then $C_n(a_1,a_2,\dots,a_k)$ is \nbcd.
\end{theorem}

Next, we present a sufficient condition for another subclass of circulant graphs to admit a \nbc.

\begin{theorem}
    Let $n$ and $s$ be multiples of $k$ and let $S = \{d_1, d_2, \dots, d_s\}$ where $1 \leq d_1 < \dots < d_s< \frac{n}{2}$.
    If exactly $\frac{s}{k}$ of the $d_t$'s are congruent to $i (\bmod\ k)$, for fixed $i\in \{1, 2,  \dots, k\}$, then the circulant graph $C_n(S)$ is \nbcd graph.
\end{theorem}
\begin{proof}
    Let $G \cong C_n(S)$ having vertex set $V(G)=\{v_1, v_2,  \dots, v_n\}$.
    Define a coloring $c \colon V(G) \to \{1, 2, \dots, k\}$ by $c(v_i) = i$, where index $i$ is taken to be modulo $k$.
    Consider a  vertex $v_r$ of $C_n(S)$.
    If $r \equiv i (\bmod\ k)$, then $r\pm d_t\equiv i\pm j (\bmod\ k)$, for all $d_t\in S$ that are congruent to $j (\bmod\ k)$.
    As $j$ takes the values from $1, 2, \dots, k$ under modulo $k$, the integer $i+j$ takes the values from $1, 2, \dots, k$ under modulo $k$. The same holds for the integer $i-j$. 
    As the number of $d_t$'s that are congruent to $j (\bmod\ k)$ is $\frac{s}{k}$, each color appears exactly $\frac{2s}{k}$ times in the neighborhood of $v_r$. 
\end{proof}

Now we give a complete characterization of \nbcd hamming graphs.
We first give a formal definition of hamming graphs.
For a set $S$, define $S^{d} := \underbrace{S \times S \times \cdots \times S}_{d-\mathrm{times}}$ to be the $d$-fold cartesian product of $S$ with itself.
A {\it hamming} graph, denoted as $H(d,k)$, is a graph with vertex set $S^{d}$, where $|S| = k$, and two vertices are adjacent if they differ in exactly one coordinate.
Note that the hamming $H(d,k)$ is $d(k-1)$-regular.
For our convenience, without loss of generality, we take $S = \{1,2, \dots, k\}$.

\begin{theorem} \label{th: hamming graphs}
The hamming graph $H(d,k)$ is a \nbcd graph if and only if $d\equiv 0\ (\bmod\ k)$.
\end{theorem}

\begin{proof}
If $d \not \equiv 0\ (\bmod\ k)$ , then the degree of any vertex of $H(d,k)$ is not a multiple of $k$ and hence by Lemma~\ref{th:prep1}, it is not a \nbcd graph.
Conversely, suppose that $d \equiv 0\ (\bmod\ k)$, so we assume $d = kn$.
We show that $H(kn,k)$ admits a \nbc by induction on $n$.
\par
For $n=1$, consider the hamming graph $H(k,k)$.
Each vertex of $H(k,k)$ is represented as a $k$-tuple $(a_1, \dots, a_k)$ where each $a_i \in S$.
For each $(a_1, a_2, \dots, a_{k-1}) \in S^{k-1}$, consider set $X^{1}_{(a_1, a_2,  \dots, a_{k-1})} = \{(a_1, a_2,  \dots, a_{k-1}, y_k) : y_k \in S\}$.
That is, each $X^1$ is a collection of vertices having  first $k-1$ coordinates the same.
There are $k^{k-1}$ such sets, each of cardinality $k$, and it gives us a partition of $S^k = V(H)$ (refer Table \ref{tab: partition of vertex set of H(4,4)}).
Again for each $(a_1, a_2, \dots, a_{k-2}) \in S^{k-2}$, let $X^{2}_{(a_1, a_2, \dots, a_{k-2})} : = \{(a_1, a_2, \dots,a_{k-2}, y_{k-1}, y_{k}) : y_{k-1}, y_k \in S\}$.
That is, each $X^2$ is a collection of vertices having  first $k-2$ coordinates the same.
There are $k^{k-2}$ such sets, each of cardinality $k^2$, and it gives a partition of the collection of all sets of the form $X^1$.
Continuing in this way, for each $r = 1, 2, \dots, k-1$, we obtain the set $X^{r}_{(a_1,a_2,\dots,a_{k-r})} = \{(a_1,a_2, \dots, a_{k-r}, y_{k-r+1}, \dots, y_k) : y_{k-r+1}, \dots, y_k \in S\}$ having  first $k-r$ coordinates the same.
Note that for each $r = 1, 2, \dots, k-1$, $|X_{(a_1,a_2,\dots,a_{k-r})}^{r}| = k^r$ and it gives us a partition of all sets of the form $X^{r-1}$ (refer Table \ref{tab: partition of vertex set of H(4,4)}).
It is straightforward to verify that for each $1 \le i \le k-1$, the set $X^i$ consists of $k$ sets of type $X^{i-1}$.
Note that for each $1 \le j \le k-1$, the tuple in the suffix of the sets of the type $X^j$ is of length $k-j$.
\par
Now we give the coloring scheme.
Consider $(1, 1,\dots, 1) \in S^{k-1}$.
We color the $k$ vertices in $X^1_{(1,1, \dots, 1)}$ by the coloring $c^{1}$ so that it is a rainbow set.
Then for $i \in S \setminus \{1\}$, we color the vertices of the form $X^1_{(1,1, \dots, 1, i)}$ by $c^{1}_{i}$ obtained by $i$th cyclic shift of colors in $c^{1}$ (recall that by our convention, the first cyclic shift of $c^1$ is $c^1$ itself).
This gives a coloring of $X^2_{(1,1, \dots, 1)}$.
We call this coloring scheme $c^2$.
Again for each $i \in S \setminus \{1\} $, we color $k^2$ vertices in $X^{2}_{(1, 1, \dots, 1, i)}$ by $c^{2}_{i}$ obtained by the $i$th cyclic shift of $c^{2}$. 
Continuing in this way, for each $1 \le j \le k-1$, we color the $k^{j}$ vertices in $X^{j}_{(1, 1, \dots, 1, i)}$ by $c^{j}_{i}$ obtained by $i$th cyclic shift of $c^{j}$ and hence we get a coloring $c^{j+1}$ of $X^{j+1}_{(1, 1, \dots, 1)}$.
In this way, after obtaining coloring $c^{k-1}$ for $X^{k-1}_{1}$, for each $i \in S \setminus \{1\}$, we apply the same coloring for $X^{k-1}_{i}$.
This gives us a coloring for $S^{k} = V(H)$.
\par
Now we show that the above coloring is \nbc.
For this consider an arbitrary vertex $(b_1, b_2, \dots, b_k)$ of $G$.
By above partitions, we have $(b_1, b_2, \dots, b_k) \in X^{1}_{(b_1, b_2,  \dots, b_{k-1})}$ and since it is a rainbow set of $k$ vertices, $(b_1, b_2,\dots, b_{k})$ has $k-1$ neighbors of distinct colors in $X^{1}_{(b_1, b_2,\dots, b_{k-1})}$.
Further, the vertex $(b_1, b_2,\dots, b_{k})$ has $k-1$ neighbors, exactly one in each of the $k-1$ other $X^1$ sets, which are in the same $X^{2}_{(b_1, b_2,\dots, b_{k-2})}$.
Continuing in this way, for each $2 \le i \le k-2$, $(b_1, b_2, \dots, b_{k}) \in X^{i}_{(b_1, b_2,\dots, b_{k-i})}$ has $k-1$ neighbors, exactly one in each of the other $X^{i}$s, which are in $X^{i+1}_{b_1, b_2,\dots, b_{k-i-1}}$.
Note that, so far, the vertex $(b_1, b_2, \dots, b_k)$ has exactly $k-1$ neighbors in each of the $k-1$ color classes, different from its own color.
Lastly, the vertex $(b_1, b_2, \dots, b_k)$ has  $k-1$ neighbors, exactly one in each of the other $X^{k-1}$ sets.
All of these neighbors have the same color as $(b_1, b_2,\dots, b_k)$, since all $X^{k-1}$ sets are colored under the same coloring scheme.
This ensures that the vertex $(b_1, b_2, \dots, b_k)$ has an equal number of neighbors in each color class.
Therefore, we conclude that $H(k,k)$ is a \nbcd graph.
This completes the base case.
\par
Now suppose $n \ge 2$.
Assume that the result is true for all hamming graphs of the form $H(kj, k)$ for all $ j < n$.
Consider the hamming graph \( H(kn, k) \).
By definition, \( H(kn, k) \) contains \( k \) vertex-disjoint copies of \( H(kn - 1, k) \), each of which in turn contains \( k \) vertex-disjoint copies of \( H(kn - 2, k) \), and so on.
In particular, \( H(kn, k) \) contains \( k^k \) vertex-disjoint copies of \( H(kn - k, k) \).
By the induction hypothesis, each copy of \( H(kn - k, k) \) admits a \nbc, say \( c_{kn-k} \).
\par
Now, consider a specific copy of \( H(kn - 1, k) \) in \( H(kn, k) \).
Within this  \( H(kn - 1, k) \), select one copy of \( H(kn  - 2, k) \), then a copy of \( H(kn - 3, k) \), and continue down to a copy \( H(kn - k, k) \).
Let us assume that this copy is colored using \( c_{kn-k} \).
In the same copy of \( H(kn - k + 1, k) \), color the remaining \( k - 1 \) copies of \( H(kn - k, k) \) using the $k-1$ colorings \( c^2_{kn-k}, c^3_{kn-k}, \dots, c^{k}_{kn-k} \) obtained by cyclic shift of \( c_{kn-k} \).
We denote the resulting coloring of this \( H(kn - k + 1, k) \) as \( c_{kn-k+1} \).
Repeat this process: use cyclic shifts \( c^2_{kn-k+1}, c^3_{kn-k+1}, \dots, c^{k}_{kn-k+1} \) to color the remaining \( k - 1 \) copies of \( H(kn - k + 1, k) \), resulting in a coloring \( c_{kn-k+2} \) of  \( H(kn - k + 2, k) \).
Continue this recursive coloring until a coloring \( c_{kn-1} \) is obtained for \( H(kn  - 1, k) \).
Apply this same coloring to each of the other \( k - 1 \) copies of \( H(kn  - 1, k) \) to obtain a complete coloring of \( H(kn, k) \).
\par
Now consider a vertex \( v \in V(H(kn, k)) \).
This vertex belongs to a specific copy of \( H(kn - k+1, k) \), say $H_v$, which lies within some \( H(kn - k + 2, k) \), which in turn lies within \( H(kn - k + 3, k) \), and so on, up to some \( H(kn  - 1, k) \).
In \(H_v =  H(kn - k+1, k) \), the vertex \( v \) has exactly one neighbor in each of the other $k-1$ color classes, as all $H(kn-k,k)$ copies receive one of the colorings from \( c_{kn-k} = c^1_{kn-k}, c^2_{kn-k}, \dots, c^{k}_{kn-k} \), all of which are cyclic shifts of the neighborhood-balanced $k$-coloring $c_{kn-k}$.
Additionally, \( v \) has \( k - 1 \) neighbors, exactly one from each of the remaining $k-1$ copies of \( H(kn - k+1, k) \) within the same \( H(kn - k + 2, k) \), each with a distinct color differing by a cyclic shift.
In general, for each \( 1 \leq i \leq k - 2 \), the vertex \( v \) has \( k - 1 \) neighbors in the corresponding \( H(kn - k + i, k) \) copies within \( H(kn - k + i + 1, k) \), where again each of these neighbors receives a distinct color via cyclic rotation.
Thus, at this stage \( v \) has exactly \( k - 1 \) neighbors of every color except its own color.
Lastly, since the coloring \( c_{kn-1} \) is applied identically across all copies of \( H(kn - 1, k) \) within \( H(kn, k) \), each vertex \( v \) has \( k - 1 \) neighbors of its own color, ensuring that $N(v)$ is equally colored.
This shows that $H(kn, k)$ is \nbcd.
\end{proof}
\begin{table}[p]
    \setlength{\tabcolsep}{5pt}
    \centering
    \centering

\begin{subtable}[t]{0.48\textwidth}
  \centering
  \caption{Partition of $X_{(1)}^{3}$}
  \renewcommand{\arraystretch}{1.5}
  \resizebox{\linewidth}{!}{%
  \begin{tabular}{|c|c|c|c|c|c|c|}
  \hline
  \multirow{16}{*}{$X_{(1)}^{3}$} & \multirow{4}{*}{$X_{(1,1)}^{2}$} & $X_{(1,1,1)}^{1}$ &\cellcolor{red!40} $(1,1,1,1)$ & \cellcolor{blue!40}$(1,1,1,2)$ & \cellcolor{green!40}$(1,1,1,3)$ & \cellcolor{yellow!40}$(1,1,1,4)$ \\ \cline{3-7}
                                  &                                  & $X_{(1,1,2)}^{1}$ &\cellcolor{blue!40} $(1,1,2,1)$ &\cellcolor{green!40} $(1,1,2,2)$ &\cellcolor{yellow!40} $(1,1,2,3)$ &\cellcolor{red!40} $(1,1,2,4)$ \\ \cline{3-7}
                                  &                                  & $X_{(1,1,3)}^{1}$ & \cellcolor{green!40}$(1,1,3,1)$ & \cellcolor{yellow!40} $(1,1,3,2)$ & \cellcolor{red!40} $(1,1,3,3)$ & \cellcolor{blue!40} $(1,1,3,4)$ \\ \cline{3-7}
                                  &                                  & $X_{(1,1,4)}^{1}$ & \cellcolor{yellow!40}$(1,1,4,1)$ & \cellcolor{red!40} $(1,1,4,2)$ & \cellcolor{blue!40} $(1,1,4,3)$ & \cellcolor{green!40} $(1,1,4,4)$ \\ \cline{2-7}
                                  & \multirow{4}{*}{$X_{(1,2)}^{2}$} & $X_{(1,2,1)}^{1}$ & \cellcolor{blue!40}$(1,2,1,1)$ & \cellcolor{green!40}$(1,2,1,2)$ &\cellcolor{yellow!40} $(1,2,1,3)$ & \cellcolor{red!40}$(1,2,1,4)$ \\ \cline{3-7}
                                  &                                  & $X_{(1,2,2)}^{1}$ & \cellcolor{green!40}$(1,2,2,1)$ & \cellcolor{yellow!40} $(1,2,2,2)$ & \cellcolor{red!40} $(1,2,2,3)$ & \cellcolor{blue!40} $(1,2,2,4)$ \\ \cline{3-7}
                                  &                                  & $X_{(1,2,3)}^{1}$ & \cellcolor{yellow!40}$(1,2,3,1)$ & \cellcolor{red!40} $(1,2,3,2)$ & \cellcolor{blue!40} $(1,2,3,3)$ & \cellcolor{green!40} $(1,2,3,4)$ \\ \cline{3-7}
                                  &                                  & $X_{(1,2,4)}^{1}$ & \cellcolor{red!40}$(1,2,4,1)$ & \cellcolor{blue!40} $(1,2,4,2)$ & \cellcolor{green!40} $(1,2,4,3)$ & \cellcolor{yellow!40} $(1,2,4,4)$ \\ \cline{2-7}
                                  & \multirow{4}{*}{$X_{(1,3)}^{2}$} & $X_{(1,3,1)}^{1}$ & \cellcolor{green!40}$(1,3,1,1)$ & \cellcolor{yellow!40} $(1,3,1,2)$ & \cellcolor{red!40} $(1,3,1,3)$ & \cellcolor{blue!40}$(1,3,1,4)$ \\ \cline{3-7}
                                  &                                  & $X_{(1,3,2)}^{1}$ & \cellcolor{yellow!40} $(1,3,2,1)$ & \cellcolor{red!40}$(1,3,2,2)$ &\cellcolor{blue!40} $(1,3,2,3)$ & \cellcolor{green!40}$(1,3,2,4)$ \\ \cline{3-7}
                                  &                                  & $X_{(1,3,3)}^{1}$ & \cellcolor{red!40}$(1,3,3,1)$ & \cellcolor{blue!40}$(1,3,3,2)$ & \cellcolor{green!40}$(1,3,3,3)$ &\cellcolor{yellow!40} $(1,3,3,4)$ \\ \cline{3-7}
                                  &                                  & $X_{(1,3,4)}^{1}$ & \cellcolor{blue!40}$(1,3,4,1)$ & \cellcolor{green!40} $(1,3,4,2)$ & \cellcolor{yellow!40} $(1,3,4,3)$ & \cellcolor{red!40}$(1,3,4,4)$ \\ \cline{2-7}
                                  & \multirow{4}{*}{$X_{(1,4)}^{2}$} & $X_{(1,4,1)}^{1}$ & \cellcolor{yellow!40}$(1,4,1,1)$ &\cellcolor{red!40} $(1,4,1,2)$ &\cellcolor{blue!40} $(1,4,1,3)$ &\cellcolor{green!40} $(1,4,1,4)$ \\ \cline{3-7}
                                  &                                  & $X_{(1,4,2)}^{1}$ & \cellcolor{red!40}$(1,4,2,1)$ &\cellcolor{blue!40} $(1,4,2,2)$ &\cellcolor{green!40} $(1,4,2,3)$ &\cellcolor{yellow!40} $(1,4,2,4)$ \\ \cline{3-7}
                                  &                                  & $X_{(1,4,3)}^{1}$ & \cellcolor{blue!40} $(1,4,3,1)$ & \cellcolor{green!40} $(1,4,3,2)$ & \cellcolor{yellow!40} $(1,4,3,3)$ & \cellcolor{red!40}$(1,4,3,4)$ \\ \cline{3-7}
                                  &                                  & $X_{(1,4,4)}^{1}$ & \cellcolor{green!40}$(1,4,4,1)$ &\cellcolor{yellow!40} $(1,4,4,2)$ &\cellcolor{red!40} $(1,4,4,3)$ &\cellcolor{blue!40} $(1,4,4,4)$ \\ \hline
  \end{tabular}%
  }
\end{subtable}%
\hfill
\begin{subtable}[t]{0.48\textwidth}
  \centering
  \caption{Partition of $X_{(2)}^{3}$}
  \renewcommand{\arraystretch}{1.5}
  \resizebox{\linewidth}{!}{%
  \begin{tabular}{|c|c|c|c|c|c|c|}
  \hline
  \multirow{16}{*}{$X_{(2)}^{3}$} & \multirow{4}{*}{$X_{(2,1)}^{2}$} & $X_{(2,1,1)}^{1}$ &\cellcolor{red!40} $(2,1,1,1)$ & \cellcolor{blue!40}$(2,1,1,2)$ & \cellcolor{green!40}$(2,1,1,3)$ & \cellcolor{yellow!40}$(2,1,1,4)$ \\ \cline{3-7}
                                  &                                  & $X_{(2,1,2)}^{1}$ &\cellcolor{blue!40} $(2,1,2,1)$ &\cellcolor{green!40} $(2,1,2,2)$ &\cellcolor{yellow!40} $(2,1,2,3)$ &\cellcolor{red!40} $(2,1,2,4)$ \\ \cline{3-7}
                                  &                                  & $X_{(2,1,3)}^{1}$ & \cellcolor{green!40}$(2,1,3,1)$ & \cellcolor{yellow!40} $(2,1,3,2)$ & \cellcolor{red!40} $(2,1,3,3)$ & \cellcolor{blue!40} $(2,1,3,4)$ \\ \cline{3-7}
                                  &                                  & $X_{(2,1,4)}^{1}$ & \cellcolor{yellow!40}$(2,1,4,1)$ & \cellcolor{red!40} $(2,1,4,2)$ & \cellcolor{blue!40} $(2,1,4,3)$ & \cellcolor{green!40} $(2,1,4,4)$ \\ \cline{2-7}
                                  & \multirow{4}{*}{$X_{(2,2)}^{2}$} & $X_{(2,2,1)}^{1}$ & \cellcolor{blue!40}$(2,2,1,1)$ & \cellcolor{green!40}$(2,2,1,2)$ &\cellcolor{yellow!40} $(2,2,1,3)$ & \cellcolor{red!40}$(2,2,1,4)$ \\ \cline{3-7}
                                  &                                  & $X_{(2,2,2)}^{1}$ & \cellcolor{green!40}$(2,2,2,1)$ & \cellcolor{yellow!40} $(2,2,2,2)$ & \cellcolor{red!40} $(2,2,2,3)$ & \cellcolor{blue!40} $(2,2,2,4)$ \\ \cline{3-7}
                                  &                                  & $X_{(2,2,3)}^{1}$ & \cellcolor{yellow!40}$(2,2,3,1)$ & \cellcolor{red!40} $(2,2,3,2)$ & \cellcolor{blue!40} $(2,2,3,3)$ & \cellcolor{green!40} $(2,2,3,4)$ \\ \cline{3-7}
                                  &                                  & $X_{(2,2,4)}^{1}$ & \cellcolor{red!40}$(2,2,4,1)$ & \cellcolor{blue!40} $(2,2,4,2)$ & \cellcolor{green!40} $(2,2,4,3)$ & \cellcolor{yellow!40} $(2,2,4,4)$ \\ \cline{2-7}
                                  & \multirow{4}{*}{$X_{(2,3)}^{2}$} & $X_{(2,3,1)}^{1}$ & \cellcolor{green!40}$(2,3,1,1)$ & \cellcolor{yellow!40} $(2,3,1,2)$ & \cellcolor{red!40} $(2,3,1,3)$ & \cellcolor{blue!40}$(2,3,1,4)$ \\ \cline{3-7}
                                  &                                  & $X_{(2,3,2)}^{1}$ & \cellcolor{yellow!40} $(2,3,2,1)$ & \cellcolor{red!40}$(2,3,2,2)$ &\cellcolor{blue!40} $(2,3,2,3)$ & \cellcolor{green!40}$(2,3,2,4)$ \\ \cline{3-7}
                                  &                                  & $X_{(2,3,3)}^{1}$ & \cellcolor{red!40}$(2,3,3,1)$ & \cellcolor{blue!40}$(2,3,3,2)$ & \cellcolor{green!40}$(2,3,3,3)$ &\cellcolor{yellow!40} $(2,3,3,4)$ \\ \cline{3-7}
                                  &                                  & $X_{(2,3,4)}^{1}$ & \cellcolor{blue!40}$(2,3,4,1)$ & \cellcolor{green!40} $(2,3,4,2)$ & \cellcolor{yellow!40} $(2,3,4,3)$ & \cellcolor{red!40}$(2,3,4,4)$ \\ \cline{2-7}
                                  & \multirow{4}{*}{$X_{(2,4)}^{2}$} & $X_{(2,4,1)}^{1}$ & \cellcolor{yellow!40}$(2,4,1,1)$ &\cellcolor{red!40} $(2,4,1,2)$ &\cellcolor{blue!40} $(2,4,1,3)$ &\cellcolor{green!40} $(2,4,1,4)$ \\ \cline{3-7}
                                  &                                  & $X_{(2,4,2)}^{1}$ & \cellcolor{red!40}$(2,4,2,1)$ &\cellcolor{blue!40} $(2,4,2,2)$ &\cellcolor{green!40} $(2,4,2,3)$ &\cellcolor{yellow!40} $(2,4,2,4)$ \\ \cline{3-7}
                                  &                                  & $X_{(2,4,3)}^{1}$ & \cellcolor{blue!40} $(2,4,3,1)$ & \cellcolor{green!40} $(2,4,3,2)$ & \cellcolor{yellow!40} $(2,4,3,3)$ & \cellcolor{red!40}$(2,4,3,4)$ \\ \cline{3-7}
                                  &                                  & $X_{(2,4,4)}^{1}$ & \cellcolor{green!40}$(2,4,4,1)$ &\cellcolor{yellow!40} $(2,4,4,2)$ &\cellcolor{red!40} $(2,4,4,3)$ &\cellcolor{blue!40} $(2,4,4,4)$ \\ \hline
  \end{tabular}%
  }
\end{subtable}%
\hfill

\vspace{1em}
\begin{subtable}[t]{0.48\textwidth}
  \centering
  \caption{Partition of $X_{(3)}^{3}$}
  \renewcommand{\arraystretch}{1.5}
  \resizebox{\linewidth}{!}{%
  \begin{tabular}{|c|c|c|c|c|c|c|}
  \hline
  \multirow{16}{*}{$X_{(3)}^{3}$} & \multirow{4}{*}{$X_{(3,1)}^{2}$} & $X_{(3,1,1)}^{1}$ &\cellcolor{red!40} $(3,1,1,1)$ & \cellcolor{blue!40}$(3,1,1,2)$ & \cellcolor{green!40}$(3,1,1,3)$ & \cellcolor{yellow!40}$(3,1,1,4)$ \\ \cline{3-7}
                                  &                                  & $X_{(3,1,2)}^{1}$ &\cellcolor{blue!40} $(3,1,2,1)$ &\cellcolor{green!40} $(3,1,2,2)$ &\cellcolor{yellow!40} $(3,1,2,3)$ &\cellcolor{red!40} $(3,1,2,4)$ \\ \cline{3-7}
                                  &                                  & $X_{(3,1,3)}^{1}$ & \cellcolor{green!40}$(3,1,3,1)$ & \cellcolor{yellow!40} $(3,1,3,2)$ & \cellcolor{red!40} $(3,1,3,3)$ & \cellcolor{blue!40} $(3,1,3,4)$ \\ \cline{3-7}
                                  &                                  & $X_{(3,1,4)}^{1}$ & \cellcolor{yellow!40}$(3,1,4,1)$ & \cellcolor{red!40} $(3,1,4,2)$ & \cellcolor{blue!40} $(3,1,4,3)$ & \cellcolor{green!40} $(3,1,4,4)$ \\ \cline{2-7}
                                  & \multirow{4}{*}{$X_{(3,2)}^{2}$} & $X_{(3,2,1)}^{1}$ & \cellcolor{blue!40}$(3,2,1,1)$ & \cellcolor{green!40}$(3,2,1,2)$ &\cellcolor{yellow!40} $(3,2,1,3)$ & \cellcolor{red!40}$(3,2,1,4)$ \\ \cline{3-7}
                                  &                                  & $X_{(3,2,2)}^{1}$ & \cellcolor{green!40}$(3,2,2,1)$ & \cellcolor{yellow!40} $(3,2,2,2)$ & \cellcolor{red!40} $(3,2,2,3)$ & \cellcolor{blue!40} $(3,2,2,4)$ \\ \cline{3-7}
                                  &                                  & $X_{(3,2,3)}^{1}$ & \cellcolor{yellow!40}$(3,2,3,1)$ & \cellcolor{red!40} $(3,2,3,2)$ & \cellcolor{blue!40} $(3,2,3,3)$ & \cellcolor{green!40} $(3,2,3,4)$ \\ \cline{3-7}
                                  &                                  & $X_{(3,2,4)}^{1}$ & \cellcolor{red!40}$(3,2,4,1)$ & \cellcolor{blue!40} $(3,2,4,2)$ & \cellcolor{green!40} $(3,2,4,3)$ & \cellcolor{yellow!40} $(3,2,4,4)$ \\ \cline{2-7}
                                  & \multirow{4}{*}{$X_{(3,3)}^{2}$} & $X_{(3,3,1)}^{1}$ & \cellcolor{green!40}$(3,3,1,1)$ & \cellcolor{yellow!40} $(3,3,1,2)$ & \cellcolor{red!40} $(3,3,1,3)$ & \cellcolor{blue!40}$(3,3,1,4)$ \\ \cline{3-7}
                                  &                                  & $X_{(3,3,2)}^{1}$ & \cellcolor{yellow!40} $(3,3,2,1)$ & \cellcolor{red!40}$(3,3,2,2)$ &\cellcolor{blue!40} $(3,3,2,3)$ & \cellcolor{green!40}$(3,3,2,4)$ \\ \cline{3-7}
                                  &                                  & $X_{(3,3,3)}^{1}$ & \cellcolor{red!40}$(3,3,3,1)$ & \cellcolor{blue!40}$(3,3,3,2)$ & \cellcolor{green!40}$(3,3,3,3)$ &\cellcolor{yellow!40} $(3,3,3,4)$ \\ \cline{3-7}
                                  &                                  & $X_{(3,3,4)}^{1}$ & \cellcolor{blue!40}$(3,3,4,1)$ & \cellcolor{green!40} $(3,3,4,2)$ & \cellcolor{yellow!40} $(3,3,4,3)$ & \cellcolor{red!40}$(3,3,4,4)$ \\ \cline{2-7}
                                  & \multirow{4}{*}{$X_{(3,4)}^{2}$} & $X_{(3,4,1)}^{1}$ & \cellcolor{yellow!40}$(3,4,1,1)$ &\cellcolor{red!40} $(3,4,1,2)$ &\cellcolor{blue!40} $(3,4,1,3)$ &\cellcolor{green!40} $(3,4,1,4)$ \\ \cline{3-7}
                                  &                                  & $X_{(3,4,2)}^{1}$ & \cellcolor{red!40}$(3,4,2,1)$ &\cellcolor{blue!40} $(3,4,2,2)$ &\cellcolor{green!40} $(3,4,2,3)$ &\cellcolor{yellow!40} $(3,4,2,4)$ \\ \cline{3-7}
                                  &                                  & $X_{(3,4,3)}^{1}$ & \cellcolor{blue!40} $(3,4,3,1)$ & \cellcolor{green!40} $(3,4,3,2)$ & \cellcolor{yellow!40} $(3,4,3,3)$ & \cellcolor{red!40}$(3,4,3,4)$ \\ \cline{3-7}
                                  &                                  & $X_{(3,4,4)}^{1}$ & \cellcolor{green!40}$(3,4,4,1)$ &\cellcolor{yellow!40} $(3,4,4,2)$ &\cellcolor{red!40} $(3,4,4,3)$ &\cellcolor{blue!40} $(3,4,4,4)$ \\ \hline
  \end{tabular}%
  }
\end{subtable}%
\hfill
\begin{subtable}[t]{0.48\textwidth}
  \centering
  \caption{Partition of $X_{(1)}^{3}$}
  \renewcommand{\arraystretch}{1.5}
  \resizebox{\linewidth}{!}{%
  \begin{tabular}{|c|c|c|c|c|c|c|}
  \hline
  \multirow{16}{*}{$X_{(4)}^{3}$} & \multirow{4}{*}{$X_{(4,1)}^{2}$} & $X_{(4,1,1)}^{1}$ &\cellcolor{red!40} $(4,1,1,1)$ & \cellcolor{blue!40}$(4,1,1,2)$ & \cellcolor{green!40}$(4,1,1,3)$ & \cellcolor{yellow!40}$(4,1,1,4)$ \\ \cline{3-7}
                                  &                                  & $X_{(4,1,2)}^{1}$ &\cellcolor{blue!40} $(4,1,2,1)$ &\cellcolor{green!40} $(4,1,2,2)$ &\cellcolor{yellow!40} $(4,1,2,3)$ &\cellcolor{red!40} $(4,1,2,4)$ \\ \cline{3-7}
                                  &                                  & $X_{(4,1,3)}^{1}$ & \cellcolor{green!40}$(4,1,3,1)$ & \cellcolor{yellow!40} $(4,1,3,2)$ & \cellcolor{red!40} $(4,1,3,3)$ & \cellcolor{blue!40} $(4,1,3,4)$ \\ \cline{3-7}
                                  &                                  & $X_{(4,1,4)}^{1}$ & \cellcolor{yellow!40}$(4,1,4,1)$ & \cellcolor{red!40} $(4,1,4,2)$ & \cellcolor{blue!40} $(4,1,4,3)$ & \cellcolor{green!40} $(4,1,4,4)$ \\ \cline{2-7}
                                  & \multirow{4}{*}{$X_{(4,2)}^{2}$} & $X_{(4,2,1)}^{1}$ & \cellcolor{blue!40}$(4,2,1,1)$ & \cellcolor{green!40}$(4,2,1,2)$ &\cellcolor{yellow!40} $(4,2,1,3)$ & \cellcolor{red!40}$(4,2,1,4)$ \\ \cline{3-7}
                                  &                                  & $X_{(4,2,2)}^{1}$ & \cellcolor{green!40}$(4,2,2,1)$ & \cellcolor{yellow!40} $(4,2,2,2)$ & \cellcolor{red!40} $(4,2,2,3)$ & \cellcolor{blue!40} $(4,2,2,4)$ \\ \cline{3-7}
                                  &                                  & $X_{(4,2,3)}^{1}$ & \cellcolor{yellow!40}$(4,2,3,1)$ & \cellcolor{red!40} $(4,2,3,2)$ & \cellcolor{blue!40} $(4,2,3,3)$ & \cellcolor{green!40} $(4,2,3,4)$ \\ \cline{3-7}
                                  &                                  & $X_{(4,2,4)}^{1}$ & \cellcolor{red!40}$(4,2,4,1)$ & \cellcolor{blue!40} $(4,2,4,2)$ & \cellcolor{green!40} $(4,2,4,3)$ & \cellcolor{yellow!40} $(4,2,4,4)$ \\ \cline{2-7}
                                  & \multirow{4}{*}{$X_{(4,3)}^{2}$} & $X_{(4,3,1)}^{1}$ & \cellcolor{green!40}$(4,3,1,1)$ & \cellcolor{yellow!40} $(4,3,1,2)$ & \cellcolor{red!40} $(4,3,1,3)$ & \cellcolor{blue!40}$(4,3,1,4)$ \\ \cline{3-7}
                                  &                                  & $X_{(4,3,2)}^{1}$ & \cellcolor{yellow!40} $(4,3,2,1)$ & \cellcolor{red!40}$(4,3,2,2)$ &\cellcolor{blue!40} $(4,3,2,3)$ & \cellcolor{green!40}$(1,3,2,4)$ \\ \cline{3-7}
                                  &                                  & $X_{(1,3,3)}^{1}$ & \cellcolor{red!40}$(1,3,3,1)$ & \cellcolor{blue!40}$(1,3,3,2)$ & \cellcolor{green!40}$(1,3,3,3)$ &\cellcolor{yellow!40} $(1,3,3,4)$ \\ \cline{3-7}
                                  &                                  & $X_{(1,3,4)}^{1}$ & \cellcolor{blue!40}$(1,3,4,1)$ & \cellcolor{green!40} $(1,3,4,2)$ & \cellcolor{yellow!40} $(1,3,4,3)$ & \cellcolor{red!40}$(1,3,4,4)$ \\ \cline{2-7}
                                  & \multirow{4}{*}{$X_{(1,4)}^{2}$} & $X_{(1,4,1)}^{1}$ & \cellcolor{yellow!40}$(1,4,1,1)$ &\cellcolor{red!40} $(1,4,1,2)$ &\cellcolor{blue!40} $(1,4,1,3)$ &\cellcolor{green!40} $(1,4,1,4)$ \\ \cline{3-7}
                                  &                                  & $X_{(1,4,2)}^{1}$ & \cellcolor{red!40}$(1,4,2,1)$ &\cellcolor{blue!40} $(1,4,2,2)$ &\cellcolor{green!40} $(1,4,2,3)$ &\cellcolor{yellow!40} $(1,4,2,4)$ \\ \cline{3-7}
                                  &                                  & $X_{(1,4,3)}^{1}$ & \cellcolor{blue!40} $(1,4,3,1)$ & \cellcolor{green!40} $(1,4,3,2)$ & \cellcolor{yellow!40} $(1,4,3,3)$ & \cellcolor{red!40}$(1,4,3,4)$ \\ \cline{3-7}
                                  &                                  & $X_{(1,4,4)}^{1}$ & \cellcolor{green!40}$(1,4,4,1)$ &\cellcolor{yellow!40} $(1,4,4,2)$ &\cellcolor{red!40} $(1,4,4,3)$ &\cellcolor{blue!40} $(1,4,4,4)$ \\ \hline
  \end{tabular}%
  }
\end{subtable}%

\label{tab:four-blocks}

    \caption{This table presents the partition of the vertex set of the Hamming graph $H(k,k)$ for $k=4$, as outlined in the proof of Theorem \ref{th: hamming graphs}.
    Each subtable illustrates the partition of the set $X_{(a)}^{3}$ for each $1 \le a \le 4$, which we refer to as a first level of partition. The second column in each subtable denotes the partition of $X_{a}^{3}$ into sets $X_{(a,b)}^{2}$ for $1 \le b \le 4$, which we refer to as the second level of partition.
    Finally, the third column in each subtable represents the partition of $X^{2}_{(a,b)}$ into sets $X^1_{(a,b,c)}$ for each $1 \le c \le 4$, which we refer to as the third level of partition.
    This table also provides a neighborhood-balanced $4$-coloring of $H(4,4)$.}
    \label{tab: partition of vertex set of H(4,4)}
\end{table}

Since the Hypercube $Q_d$ is a special case of the Hamming graph when $k=2$, we have the following corollary which already appeared in \cite{collin-cnbc}.

\begin{corollary}
   A hypercube $Q_d$ admits a neighborhood-balanced $2$-coloring if and only if $d$ is even.
\end{corollary}


\subsection{\small Neighborhood-balanced \texorpdfstring{k}{$k$}-colored graphs having unequal sizes of color classes}\label{sec:unequalcolorclass}

We know that for \nbcd regular graphs the color classes are of the same size (see Corollary \ref{th: regular graphs necessary condition}).
This need not be the case in general.
In this section, we give a way to start with a \nbcd graph and construct a new \nbcd graph that has fewer vertices of one color compared to the other colors.
Similar constructions for $k=2$ are shown in~\cite{collin-cnbc}.

\begin{definition}[\(2k-1\)-vertex addition]\label{def:vertexaddition}
Let $G$ be a \nbcd graph without isolated vertices and $\{u_i, v_i : 1 \le i \le k\} \subseteq V(G)$ be such that the color of $u_i$ is the same color as the color of $v_i$.
A {\it $(2k-1)$-vertex addition} at $\{u_1,u_2, \dots,u_{k}, v_1, v_2, \dots, v_{k}\}$ is the operation of adding vertices $w, a_1, b_1, a_2, b_2, \dots, a_{k-1}, b_{k-1}$ such that $w$ is adjacent to all of $u_1, u_2, \dots, u_{k}, v_1, v_2,  \dots, v_{k}$ and $a_i$ is adjacent to $u_1,u_2,  \dots, u_{k}$ and $b_i$ is adjacent to $v_1, v_2, \dots, v_{k},$ for $1 \leq i \leq k-1$ and assigning the color $k$ to the vertex $w$ and color $i$  to vertices $a_i$ and $b_i$, for $1 \leq i \leq k-1$.
\end{definition}

A $5$-vertex addition at the neighborhood balanced $3$-colored graph $K_{3,3}$ is shown in Figure \ref{fig:vertexadd}.
\begin{figure}[!htbp]
    \centering
    \begin{tikzpicture}

  \node[circle, draw=black, thick, minimum size=3mm, label={[xshift=0pt, yshift=-24pt]{$u_1$}}, inner sep=0pt, fill=red]   (u1) at (0, 2)  {};
  \node[circle, draw=black, thick, minimum size=3mm, label={[xshift=0pt, yshift=-24pt]{$u_2$}}, inner sep=0pt, fill=blue]  (u2) at (0, 0)  {};
  \node[circle, draw=black, thick, minimum size=3mm, label={[xshift=0pt, yshift=-24pt]{$u_3$}}, inner sep=0pt, fill=green] (u3) at (0,-2)  {};

  \node[circle, draw=black, thick, minimum size=3mm, label={[xshift=-12pt, yshift=-12pt]{$a_1$}}, inner sep=0pt, fill=red] (a1) at (-3,1)  {};
  \node[circle, draw=black, thick, minimum size=3mm, label={[xshift=-12pt, yshift=-12pt]{$a_2$}}, inner sep=0pt, fill=blue] (a2) at (-3,-1)  {};

  \node[circle, draw=black, thick, minimum size=3mm, label={[xshift=0pt, yshift=-24pt]{$v_1$}}, inner sep=0pt, fill=red] (v1) at (4, 2)  {};
  \node[circle, draw=black, thick, minimum size=3mm, label={[xshift=0pt, yshift=-24pt]{$v_2$}}, inner sep=0pt, fill=blue] (v2) at (4, 0)  {};
  \node[circle, draw=black, thick, minimum size=3mm, label={[xshift=0pt, yshift=-24pt]{$v_3$}}, inner sep=0pt, fill=green] (v3) at (4,-2) {};

   \node[circle, draw=black, thick, minimum size=3mm, label={[xshift=12pt, yshift=-12pt]{$b_1$}}, inner sep=0pt, fill=red] (b1) at (7,1)  {};
  \node[circle, draw=black, thick, minimum size=3mm, label={[xshift=12pt, yshift=-12pt]{$b_2$}}, inner sep=0pt, fill=blue] (b2) at (7,-1)  {};

  \node[circle, draw=black, thick, minimum size=3mm, label={[xshift=12pt, yshift=-12pt]{$w$}}, inner sep=0pt, fill=green] (w) at (2,4)  {};

  \draw[thick] (u1) -- (v1);
  \draw[thick] (u1) -- (v2);
  \draw[thick] (u1) -- (v3);

  \draw[thick] (u2) -- (v1);
  \draw[thick] (u2) -- (v2);
  \draw[thick] (u2) -- (v3);

  \draw[thick] (u3) -- (v1);
  \draw[thick] (u3) -- (v2);
  \draw[thick] (u3) -- (v3);

   \draw[thick] (a1) -- (u1);
  \draw[thick] (a1) -- (u2);
  \draw[thick] (a1) -- (u3);

   \draw[thick] (a2) -- (u1);
  \draw[thick] (a2) -- (u2);
  \draw[thick] (a2) -- (u3);

   \draw[thick] (b1) -- (v1);
  \draw[thick] (b1) -- (v2);
  \draw[thick] (b1) -- (v3);

   \draw[thick] (b2) -- (v1);
  \draw[thick] (b2) -- (v2);
  \draw[thick] (b2) -- (v3);

  \draw[thick] (w) -- (v1);
  \draw[thick] (w) -- (v2);
  \draw[thick] (w) -- (v3);
  \draw[thick] (w) -- (u1);
  \draw[thick] (w) -- (u2);
  \draw[thick] (w) -- (u3);

\end{tikzpicture}
    \caption{A $5$-vertex addition at the neighborhood balanced $3$-colored graph $K_{3,3}$ as described in Definition \ref{def:vertexaddition}. }
    \label{fig:vertexadd}
\end{figure}
It is straightforward to verify that a graph obtained by $(2k-1)$-vertex addition as defined above from a \nbcd is \nbcd.
Also, the $(2k-1)$-vertex addition adds one additional vertex of one color and two additional vertices, each of the other $(k-1)$ colors.
Hence, we have the following proposition.
\begin{proposition}\label{vertexaddition}
    Given a \nbc of a graph $G$, if $G'$ is the graph resulting from $(2k-1)$-vertex addition at a set $\{u_1,u_2,\dots,u_{k},v_1,v_2,\dots, v_{k}\}$, then $G'$ is a \nbcd graph. Moreover, $G'$ has one additional vertex of one color and two additional vertices, each of the other $(k-1)$ colors.
\end{proposition}

\begin{corollary}
    There exist \nbcd graphs that have a \nbc with an arbitrary fewer vertices of one color than the other colors.
    Moreover, every \nbcd graph is an induced subgraph of such a graph. 
\end{corollary}
\begin{proof}
    Let $G$ be a \nbcd graph and fix a \nbc of $G$.
    By Proposition \ref{vertexaddition}, a graph $G'$ obtained from $G$ by $(2k-1)$-vertex addition is \nbcd graph and has one additional vertex of only one color, say $1$, and two additional vertices of each of the other $k-1$ colors.
    We again do $(2k-1)$-vertex addition at $G'$ to obtain another graph $G''$ that now has two additional vertices (than $G$) of color $1$, and four additional vertices (than $G$) of other $k-1$ colors.
    Repeating such $(2k-1)$-vertex additions, finally, we obtain a \nbcd graph that has arbitrarily fewer vertices of color $1$ as compared to vertices of other $k-1$ colors.
\end{proof}

 \section{Hardness Results} \label{sec:hardness}





In this section, we show that \probNBCshort is \NP-complete by giving a reduction from \textsc{$k$-PC}.

\begin{figure}[!h]
    \centering
   \includegraphics[width=0.27\textwidth]{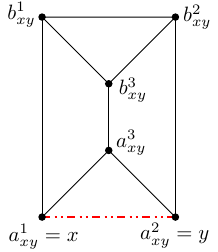}
    \caption{Gadget to convert 3-PC to 3-NBC.}
    \label{fig:3nbcgadget}
\end{figure}

\begin{figure}[h]
    \centering
   \includegraphics[width=0.27\textwidth]{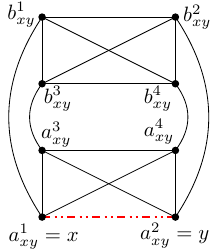}
    \caption{Gadget to convert 4-PC to 4-NBC}    \label{fig:4nbcgadget}
\end{figure}

Our reduction is based on the following gadget.
Given two vertices $x$ and $y$, we define the graph $H_{xy}$ as follows.
Let $V(H_{xy})=A_{xy}\cup B_{xy}$, where
\begin{equation*}
    A_{xy}=\{a^{1}_{xy}=x,a^{2}_{xy}=y,a^{3}_{xy},\ldots,a^{k}_{xy}\}, \quad B_{xy}=\{b^{1}_{xy},b^{2}_{xy},\ldots,b^{k}_{xy}\}
\end{equation*}
The edge set  of $H_{xy}$ is defined by
\begin{equation*}
E(H_{xy})=(A_{xy}\times A_{xy})\cup(B_{xy}\times B_{xy})\cup
\{a^{i}_{xy} b^{i}_{xy}:1\le i\le k\} - \{xy\}.
\end{equation*}
In other words, both $A_{xy}$ and $B_{xy}$ induce cliques of size  $k$, and for each $i\in[k]$, the vertex $a_{xy}^{i}$ is adjacent to the corresponding vertex $b_{xy}^{i}$. The only missing edge in the clique induced by $A_{xy}$ is the edge $xy$. 
Observe that, every vertex of $H_{xy}$, except $x$ and $y$, has degree exactly $k$ (see Figures \ref{fig:3nbcgadget} and \ref{fig:4nbcgadget} ).

\begin{figure}[h]
    \centering
   \includegraphics[width=0.36\textwidth]{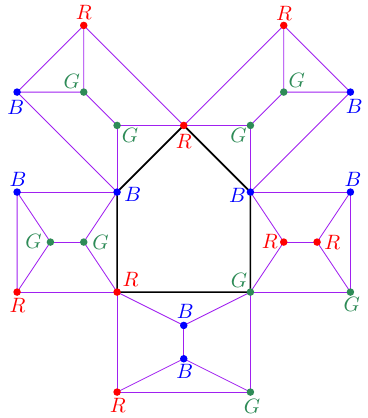}
    \caption{The 3-NBC graph $C_5'$ constructed from 3-PC graph $C_5$.}
    \label{fig:3NBC_graphs}
\end{figure}

Given an instance $(G,k)$ of \textsc{$k$-PC}, we construct a graph $G'$ as follows.
For every edge $xy \in E(G)$, attach a copy of the gadget $H_{xy}$ by identifying the vertices $x, y$ of the gadget with the vertices $x$ and $y$ of $G$.
Thus,
$$
V(G') = \bigcup_{xy \in E(G)} V(H_{xy}), \quad \text{and} \quad E(G')=E(G)\cup\left(\bigcup_{xy \in E(G)} E(H_{xy}) \right).
$$
Refer Figure \ref{fig:3NBC_graphs} for construction of 3-NBC graph $C_5'$ from 3-PC graph $C_5$.

\begin{theorem}\label{the:knbcHardness}
\probNBC problem is \NP-complete for every integer $k\ge3$.
\end{theorem}

\begin{proof}
Throughout the proof we assume that $k \geq 3$.
We prove the result by a polynomial-time reduction from the \textsc{$k$-PC} problem.
Let $(G,k)$ be an instance of \textsc{$k$-PC}, and construct the graph $G'$ as described above.
Since one gadget $H_{xy}$ of constant size (depending only on $k$) is added for every edge $xy$ of $G$, the graph $G'$ can clearly be constructed in polynomial time.
Further, given a coloring of a graph, one can verify in polynomial time whether every vertex has an equal number of neighbors of each color.
Hence, the \probNBCshort problem belongs to the class \NP.
We now show that $(G,k)$ is a \yes-instance of \textsc{$k$-PC} if and only if $(G',k)$ is a \yes-instance of \probNBCshort.

\medskip

\noindent\textbf{($\Rightarrow$)} Suppose that $(G,k)$ is a \textsc{Yes}-instance of \textsc{$k$-PC}. Let $c:V(G)\rightarrow[k]$
be a $k$-proper coloring of $G$.
We extend $c$ to a neighborhood balanced $k$-coloring $c'$ of $G'$ as follows (see Figure \ref{fig:3NBC_graphs}).

\begin{enumerate}[noitemsep]
\item Every vertex of $G$ retains the color assigned by $c$.

\item For every gadget $H_{xy}$, assign the vertices
$a^{3}_{xy},\ldots,a^{k}_{xy}$ the remaining $k-2$ colors that are different from those assigned to $x$ and $y$.

\item For each $1\le i\le k$, assign $b^{i}_{xy}$ the same color as $a^{i}_{xy}$.
\end{enumerate}

We show that $c'$ is a neighborhood balanced $k$-coloring of $G'$.
First , consider a vertex $b^{i}_{xy}$.
Its neighborhood consist of the vertex $a^{i}_{xy}$ together with the remaining $k-1$ vertices of the clique $B_{xy}$. Since each color appears exactly once among these vertices, the neighborhood of $b^{i}_{xy}$ contains exactly one vertex of each color.
Next, consider a vertex $a^{i}_{xy}$, where $3\le i\le k$. Its neighborhood consists of the vertices  $x$, $y$, $b^{i}_{xy}$ and $a^{j}_{xy}$ for all $j\neq i$ with $3\le j\le k$.
Since $c$ is a proper coloring of $G$, the vertices $x$ and $y$ receive distinct colors under $c$ in $G$.
Hence, by definition of $c'$, the vertices $x$ and $y$ receives distinct colors under $c'$ in $G'$.
Moreover, the vertices $a^{3}_{xy},\ldots ,a^{k}_{xy}$ are assigned the remaining $k-2$ colors, while $b^{i}_{xy}$ receives the same color as $a^{i}_{xy}$.
Consequently, every color appears exactly once in the neighborhood of $a^{i}_{xy}$.

Finally, consider a vertex $x \in V(G)$.
Let $y$ be an arbitrary neighbor of $x$ in $G$. Within the gadget $H_{xy}$, the neighbors of $x$ are the vertices $y$, $a^{3}_{xy},\ldots,a^{k}_{xy}$, and $b^{1}_{xy}$.
These vertices contribute exactly one occurrence of each of the $k$ colors: the vertex $y$ contributes its own color, the vertices $a^{3}_{xy},\ldots,a^{k}_{xy}$ contribute the remaining $k-2$ colors, and $b^{1}_{xy}$ contributes the color of $x$. Since this is true for every neighbor $y$ of $x$, every incident gadget contributes exactly one vertex of each color to the neighborhood of $x$ in $G'$.
Hence, the neighborhood of $x$ in $G'$ contains the equal number of vertices of each color.
Thus, we have shown that every vertex in $G'$ has equal number of neighbors of each color.
Therefore, $c'$ is a neighborhood balanced $k$-coloring of $G'$, and so $(G',k)$ is a \yes-instance of \probNBCshort.

\medskip

\noindent\textbf{($\Leftarrow$)} Suppose that $(G',k)$ is a \yes-instance of \probNBCshort. Let $c':V(G') \rightarrow [k]$
be a \nbc of $G'$.
Consider any gadget $H_{xy}$ for $xy \in E(G)$.
For every $3\le i\le k$, the degree of vertex $a^{i}_{xy}$ in $G'$ is $k$.
Since $c'$ is a neighborhood-balanced coloring, each of its $k$ neighbors must receive a distinct color.
As both $x=a^{1}_{xy}$ and $y=a^{2}_{xy}$ are neighbors of $a^{i}_{xy}$, they must receive distinct colors.

Since the above argument holds for every edge $xy \in E(G)$, every pair of adjacent vertices of $G$ receives different colors under the restriction of $c'$ to $V(G)$.
Hence, the restriction of $c'$ to $V(G)$ is a $k$-proper coloring of $G$.
Therefore, $(G,k)$ is a \yes-instance of \textsc{$k$-PC}.
This completes the proof.
\end{proof}

\section{Conclusion}

In this article, we defined the concept of neighborhood balanced $k$-coloring of graphs, which is a generalization of neighborhood balanced $2$-coloring of graphs introduced by Freyberg and Marr \cite{nbc}. Initially, we gave some characteristics of graphs that admit such a labeling. Then we showed how more \nbcd graphs can be constructed from the existing \nbcd using various graph operations. Further, we presented several \nbcd graphs.
We now pose the following problems related to our work.
\begin{enumerate}
    \item Characterize the regular graphs that admit a \nbc.
    \item Given a \nbcd graph $G$ and an integer $n$, satisfying conditions as in Theorem \ref{th: unionoversepset}, characterize the dependent set $S$ (that is, the subgraph $G\left< S \right>$ induced by $S$) such that the union of $G$ over $S$ is a \nbcd graph.
\end{enumerate}

\vspace{1cm}
\noindent \textbf{\large Fundings:}\\
 The first author is thankful to the Goa State Research Foundation (GSRF), Goa, India, for the financial support to carry out this work under the minor research project MIN2025035.\\
The second and fourth Authors are thankful to the National Board for Higher Mathematics(NBHM), Govt. of India, for the financial support to carry out research under the project Ref. No.: 02011/21/2025/NBHM-RP/RD-II/9821.\\

\noindent\textbf{\large Data Availability:}\\
No databases were used, and no data collection was done in this work.\\

\noindent\textbf{\large Conflict of Interest:}\\
The authors declare that they have no conflict of interest.\\

\noindent\textbf{\large Author Contributions:}\\
Each author contributed ideas and worked together to obtain the results. The first draft of the manuscript was written by Maurice G. Almeida, who is the corresponding author. All authors read and approved the final manuscript.


\bibliographystyle{plainurl}
\bibliography{bibliography.bib}

\end{document}